\documentclass[10pt,a4paper]{article}
\usepackage[utf8]{inputenc}
\usepackage[T1]{fontenc}
\usepackage[russian]{babel}
\usepackage{amsmath, amsfonts, amssymb, graphicx}
\usepackage[left=1.50cm, right=1.50cm, top=2.00cm, bottom=1.30cm]{geometry}
\usepackage{mathrsfs, amsthm}
\usepackage{verbatim}
\usepackage{fancyhdr}
\newtheorem{theorem}{Теорема}[section]

\theoremstyle{definition}
\newtheorem{definition}[theorem]{Определение}
\newtheorem{example}[theorem]{Пример}
\newtheorem{algorithm}[theorem]{Алгоритм}
\newtheorem{hypothesis}[theorem]{Гипотеза}

\theoremstyle{remark}
\newtheorem{remark}[theorem]{Замечание}
\newtheorem{assertion}[theorem]{Утверждение}
\newtheorem{corollary}[theorem]{Заключение}

\providecommand{\keywords}[1]
{
	\small	
	\textbf{Ключевые слова: } #1
}

\providecommand{\udk}[1]
{
	\small	
	\textbf{УДК } #1
}

\title{Realization of Integrable Hamiltonian Systems by Billiard Books
	\footnotetext{
		Раздел \ref{s:2:1:graphs} выполнен при поддержке гранта РНФ ($N$ 17-11-01303) в МГУ имени М.В.Ломоносова, раздел \ref{s:3:local} частично выполнен при поддержке гранта РФФИ ($N$ 19-01-00775 А) в МГУ имени М.В.Ломоносова, раздел \ref{s:4:exotic} выполнен при поддержке гранта РНФ ($N$ 20-71-00155) в МГУ имени М.В. Ломоносова, раздел \ref{s:3:local} частично выполнен в Московском центре фундаментальной и прикладной математики.
}}
\author{
	A.\,T.~Fomenko\footnote{atfomenko@mail.ru, Chair of Differential Geometry and Applications, Faculty of Mechanics and Mathematics, Lomonosov Moscow State University.},
	I.\,S.~Kharcheva\footnote{irina.harcheva1@yandex.ru, Chair of Differential Geometry and Applications, Faculty of Mechanics and Mathematics, Lomonosov Moscow State University.},
	V.\,A.~Kibkalo\footnote{slava.kibkalo@gmail.com, Chair of Differential Geometry and Applications, Faculty of Mechanics and Mathematics, Lomonosov Moscow State University; Moscow Center for Fundamental and Applied Mathematics.}
}
\date{}

\begin{document}
	\maketitle
	
	\begin{abstract}
		В работе обсуждается гипотеза А.~Т.~Фоменко о реализации интегрируемыми биллиардами топологии слоений Лиувилля гладких и аналитических интегрируемых гамильтоновых систем. Изложен алгоритм Ве\-дюш\-киной--Харчевой реализации биллиардными книжками 3-атомов, значительно упростившийся благодаря его формулировке в терминах $f$-гра\-фов. Отметим, что произвольный тип базы слоения Лиувилля на всей изоэнергетической 3-поверхности был также реализован В.~В.~Ведюшкиной и И.~С.~Харчевой с помощью другого алгоритма. В работе он наглядно проиллюстрирован на примере реализации инварианта известной интегрируемой системы Жуковского (случай Эйлера с гиростатом) в определенной зоне энергии. Как оказалось, при этом было реализовано и само слоение Лиувилля, а не только класс его базы, т.~е. получена лиувиллева эквивалентность биллиардной и механической систем. Затем изложены результаты В.~В.~Ведюшкиной и В.~А.~Кибкало по построению биллиардов с произвольными числовыми инвариантами. Далее для биллиардных книжек без потенциала с определенным свойством доказано существование инварианта Фоменко--Цишанга и их принадлежность к классу топологически устойчивых систем. В конце приведен пример, когда добавление потенциала Гука к плоскому биллиарду порождает расщепляющуюся невырожденную 4-особенность ранга 1.
		
		
		\udk{517.938.5.}

	\end{abstract}
	
	\section{Введение}
	
	Новые результаты, полученные при исследовании вопроса об интегрируемости и неинтегрируемости биллиардных систем, в значительной мере прояснили связь между этим явлением и устройством границы биллиардного стола.
	
	\begin{remark} Интегрируемость биллиарда может рассматриваться в различных смыслах (как пример, локальная, полиномиальная, алгебраическая). В настоящей работе под интегрируемыми биллиардами мы будем понимать полиномиально интегрируемые биллиарды: системы с дополнительным интегралом, полиномиальным по компонентам импульса.
	\end{remark}
	
	Как оказалось, для нескольких важных классов биллиардных систем верна гипотеза Биркгофа, причем как в смысле полиномиальной интегрируемости  (см. работу А.~А.~Глуцюка \cite{Gluts}, также работы С.~В.~Болотина \cite{Bolot}, А.~Е.~Миронова и М.~Бялого \cite{MirByal1}, \cite{MirByal2}), так и ее локальная версия  (работы В.~Ю.~Калошина, С.~Соррентино и др. \cite{Sorrent1}, \cite{Sorrent2}). Иными словами, была подтверждена глубокая связь интегрируемости биллиарда на односвязном пространстве постоянной кривизны с классом софокусных квадрик или концентрических окружностей.
	
	Софокусные квадрики на плоскости образуют семейство \eqref{eq_confocal} c параметром $0 \le \lambda \le a$, где $0 < b < a$ --- полуоси эллипса $\lambda = 0$. Вдоль траекторий биллиарда в области, ограниченной их дугами, сохраняется интеграл $\Lambda$ \eqref{eq_integral}. Его значением в точке $(x, y, v_x, v_y) \in T^{*} \mathbb{R}^2$ является параметр $\lambda$ той квадрики семейства, которой касаются все звенья данной траектории.
	\begin{equation}\label{eq_confocal}
		(b - \lambda) x^2 + (a - \lambda) y^2 = (a - \lambda) (b - \lambda).
	\end{equation}
	\begin{equation}\label{eq_integral}
		\Lambda = \frac{-(x v_y - y v_x)^2 + v_x^2 b + v_y^2 a}{v_x^2 + v_y^2}.
	\end{equation}
	С другой стороны, В.~В.~Ведюшкиной удалось существенно расширить класс интегрируемых биллиардов, определив биллиардные системы на клеточных комплексах специального вида, склеенных из плоских софокусных интегрируемых биллиардов по их одинаковым дугам границы. Общая конструкция \textit{биллиардных книжек} --- таких столов, и динамических систем на них --- была введена В.~В.~Ведюшкиной в~\cite{VedKha}, а сама идея склейки более узкого класса \textit{топологических биллиардов} ---  в работе \cite{FokDomain}). При этом, фазовое пространство и неособый трехмерный уровень постоянной энергии (называемый изоэнергетической поверхностью $ Q^3 $) биллиардных книжек остаются кусочно-гладкими многообразиями, как и в случае плоских софокусных биллиардов (см.~\cite{Kha}).
	
	\begin{remark} Отметим, что комбинация склейки комплекса из плоских столов и введения одинаковых подходящих потенциалов на них (например, потенциала Гука, если границы --- софокусные квадрики) или одинакового постоянного по времени и пространству магнитного поля (для круговых биллиардов) сохраняет интегрируемость системы.
	\end{remark}
	
	Получаемый комплекс (биллиардная книжка) уже не обязан быть изометрически вложимым целиком в плоскость $Oxy$ с евклидовой метрикой. Тем не менее, определена его проекция на $Oxy$. Она является изометрией на каждой открытой 2-клетке комплекса, т.~е. на \textit{внутренности} каждого из склеиваемых столов элементарных биллиардов. Проекция определяет закон движения материальной частицы (биллиардного шара) на открытых 2-клетках (внутри листов книжки) комплекса благодаря закону движения на плоскости.
	
	Ребра комплекса (его 1-клетки) проецируются биективно на гладкие дуги кусочно-гладкой границы 2-клетки (элементарного биллиарда). Каждому ребру приписана перестановка на склеиваемых по нему столах. Она задает переход материальной частицы с листа на лист после достижения данного ребра (называемого также \textit{корешком книжки}), т.~е. после удара о границу. Значение интеграла $\Lambda$ не меняется при переходе шара с листа на лист.
	
	Оказалось, что склеивая представителей весьма узкого класса (конечное число типов плоских интегрируемых биллиардов, ограниченных дугами софокусных квадрик), можно реализовать топологические инварианты (построенные в \cite{FomDAN86}-\cite{FZ91}, подробнее см. \cite{BolFom}) слоений Лиувилля достаточно широкого класса  интегрируемых гамильтоновых систем с двумя степенями свободы (далее ИГС) \cite{FokFomDAN}-\cite{saddle21}, а также их особенностей (см. \cite{VedKha} и \cite{FomVedKha}). Вычислению инвариантов систем динамики и физики посвящены, например, работы \cite{BFR}-\cite{Kib20}.
	
	\subsection{Гипотеза Фоменко о биллиардах}
	
	В работе \cite{VedFom19} В.~В.~Ведюшкиной и А.~Т.~Фоменко сделан обзор полученных результатов и открытых задач. Там же, на основе полученных ранее результатов, А.~Т.~Фоменко сформулировал фундаментальную гипотезу о реализации (моделировании) ИГС интегрируемыми биллиардовами. Под реализацией некоторых инвариантов ИГС мы понимаем построение биллиарда, слоение Лиувилля которого имеет такое же значение данного инварианта.
	
	\begin{hypothesis}[А.\,Т.~Фоменко] \label{global_hypothesis}
		Подходящим классом интегрируемых биллиардов можно реализовать:
		\begin{itemize}
			\item
			(гипотеза $\textbf{А}$) любой атом или, другими словами, любую типичную невырожденную бифуркацию ранга 1 двумерных торов Лиувилля;
			\item
			(гипотеза $\textbf{B}$) любую грубую молекулу, или, другими словами, базу любого слоения Лиувилля на 3-мно\-го\-обра\-зии;
			\item
			(гипотеза $\textbf{C}$) любую меченую молекулу (инвариант Фоменко-Цишанга), или, другими словами, любое слоение Лиувилля  на 3-многообразии;
			\item
			(гипотеза $\textbf{D}$) любое замкнутое неособое изоэнергетическое 3-многообразие любой невырожденной интегрируемой гамильтоновой системы (с точностью до гомеоморфности).
		\end{itemize}
	\end{hypothesis}
	
	Также в \cite{VedFom19} были выделены классы биллиардов I-VIII (элементарные и топологические биллиарды, биллиардные книжки, биллиарды на плоскости Минковского, биллиарды с магнитным полем, биллиарды с потенциалом, объединяющий их класс биллиардов), в которых был начат поиск подходящих систем.
	
	Раздел $\mathbf{C}$ гипотезы Фоменко является наиболее сильным. Он предполагает возможность по произвольному выбранному инварианту Фоменко-Цишанга построить биллиард, чье слоение Лиувилля имеет такой же инвариант.
	
	Напомним, что инвариант Фоменко-Цишанга классифицирует интегрируемые системы с двумя степенями свободы с точностью до лиувиллевой эквивалентности. А именно, с точностью до послойного гомеоморфизма слоений Лиувилля систем в выбранных зонах энергии, с дополнительным условием сохранения ориентации критических окружностей, задаваемой на них гамильтоновым векторным полем. Инвариант представляет собой граф, вершины которого соответствуют бифуркациям ранга 1 и снабженный некоторыми числовыми метками $ r, \varepsilon, n $. Такой граф называется меченой молекулой.
	
	Если гипотеза $ \mathbf{C}$ верна, то множество классов лиувиллевой эквивалентности (гладких или вещественно-аналитических) ИГС с двумя степенями свободы, имеющих боттовские особенности ранга 1, являющихся топологически устойчивыми и нерезонансными \cite{BolFom}, вкладывается в множество классов кусочно-гладкой лиувиллевой эквивалентности интегрируемых биллиардов. Иными словами, в этом случае интегрируемые биллиарды \textit{моделируют} все множество указанных интегрируемых систем, или его существенную часть.
	
	Другие разделы $\mathbf{A}$, $ \mathbf{B}$, $ \mathbf{D}$ гипотезы имеют собственный геометрический смысл, являясь при этом необходимыми условиями для истинности раздела $ \mathbf{C}$. Другим таким условием является \textit{локальная гипотеза Фоменко}, сформулированная А.~Т.~Фоменко позднее в \cite{DAN20}.
	
	Разделы $ \mathbf{B}$ и $ \mathbf{D}$ гипотезы Фоменко предполагают реализацию ИГС биллиардами с точностью до \textit{более слабых} эквивалентностей, чем лиувиллева эквивалентность в разделе $ \mathbf{C}$.
	
	Для раздела $ \mathbf{B}$ имеется в виду реализация биллиардом слоения Лиувилля с произвольной базой (задаваемой \textit{грубой молекулой}, т.~е. описанным выше графом, лишенным \textit{числовых} меток, он же называется \textit{инвариантом Фоменко}).
	
	Раздел $ \mathbf{D}$ предполагает возможность реализовать (как трехмерное инвариантное многообразие биллиарда) произвольное трехмерное многообразие класса $(H)$. Это есть класс изоэнергетических многообразий для гладких невырожденных ИГС с двумя степенями свободы, совпадающий (как показано А.~Т.~Фоменко, Х.~Цишангом в \cite{FomDAN86} и \cite{FZ87}) с известным классом граф-многообразий Вальдхаузена, см.~\cite{Wald1}, \cite{Wald2}.
	
	В то же время раздел $ \mathbf{A}$ и локальная гипотеза предполагают реализуемость произвольных \textit{составных частей} инварианта Фоменко-Цишанга.
	
	Другой интересный вопрос, поставленный в \cite{VedFom19} --- о возможности реализации (в смысле послойной гомеоморфности) биллиардами \textit{четырехмерных} особенностей ранга 0 ИГС на всём симплектическом многообразии (а не только слоение Лиувилля на 3-подмногообразии). Как показано В.А. Кибкало в \cite{rankzero}, слоение подходящего биллиарда (биллиардной книжки) с отталкивающим потенциалом Гука позволяет реализовать как примеры особенностей центр-центр, центр-седло и седло-седло, так и произвольный класс особенностей типа центр-седло. Согласно подсчету С.В.Пустовойтова (об этом см. \cite{MagPot}), переход от софокусного биллиарда к круговому биллиарду с потенциалом позволяет реализовать особенность типа фокус-фокус.
	
	Данная работа продолжается, и есть основания полагать, что любую невырожденную особенность ранга 0 ИГС с 2 степенями свободы можно будет реализовать биллиардом с потенциалом. В системах с большим числом степеней свободы возникают особенности из других ``смешанных'' серий, например, седло-фокус. С точностью до полулокальной эквивалентности (послойной гомеоморфности) они были недавно классифицированы И.~К.~Козловым и А.~А.~Ошемковым \cite{KozOsh}.
	
	\subsection{Содержание работы}
	\subsubsection{Раздел \ref{s:2:1:graphs}} В доказанном В.~В.~Ведюшкиной и И.~С.~Харчевой~\cite{VedKha} разделе $ \mathbf{A}$ гипотезы Фоменко утверждается реализация произвольного  3-атома. Напомним, что так называют произвольную трехмерную боттовскую особенность ранга 1. Cимвол 3-атома приписывают каждой вершине графа Риба слоения Лиувилля на $Q^3$, получая инвариант Фоменко (или грубую молекулу) данного слоения --- классифицирующий инвариант его базы.
	
	Для реализации 3-атомов В.~В.~Ведюшкиной был введен~\cite{VedKha} класс биллиардных книжек --- биллиардных систем на клеточных комплексах, ребрам которого (1-клеткам) приписаны перестановки на склеиваемых столах (2-клетках). Подробности определения напомним ниже.
	
	Алгоритм Ведюшкиной-Харчевой формулируется в алгебраических терминах перестановок. Вместе с тем, удобной техникой для работы с атомами является подход $f$-графов, предложенный А.~А.~Ошемковым \cite{Oshemkov}. В частности, он позволяет эффективно перечислять и сравнивать особенности. Другим известным подходом является топологическое описание 3-атомов в терминах лент и крестов в их базе или дубле (для атомов со звездочками). В разделе \ref{s:2:1:graphs} мы даем представление результата работы~\cite{VedKha} в существенно более наглядных терминах $f$-графов, иллюстрируя его на примере 3-атома $B^*$.
	
	\subsubsection{Раздел \ref{s:2:molec}}  Гипотеза $ \mathbf{B}$ также верна. В работе \cite{VedKha2} В.~В.~Ведюшкиной и И.~С.~Харчевой был предложен  алгоритм построения биллиардной книжки по произвольной выбранной грубой молекуле (инвариантом базы слоения Лиувилля). Полученный алгоритм обсуждается в разделе \ref{s:2:molec} и применяется к реализации известной системы Жуковского из динамики твердого тела. Оказалось, что построенная по алгоритму биллиардная книжка реализует не только базу, но и само слоение Лиувилля, т.~е. две разные по природе системы лиувиллево эквивалентны друг другу.
	
	\subsubsection{Раздел \ref{s:3:local}} В локальной гипотезе Фоменко ставится следующий вопрос: какие значения  числовых меток и типы слоения Лиувилля ``вблизи'' атома, семьи или ребра молекулы могут быть реализованы биллиардами. Как оказалось, произвольное значение метки $n$ реализуется на некоторой семье седловых атомов некоторой молекулы \cite{DAN20, VedKib20}, а произвольное значение пары меток $(r, \varepsilon)$ реализуется на некотором ребре некоторой молекулы. Напомним, что семьей молекулы называют связный подграф молекулы, прообраз которого имеет структуру единого расслоения Зейферта, а все атомы являются седловыми (гиперболическими).
	
	\subsubsection{Раздел \ref{s:4:exotic}} В этом разделе мы обсудим, какими общими свойствами могут обладать слоения и особенности интегрируемых биллиардов.
	
	Во-первых, системы биллиардов на клеточных комплексах и их модификации являются, вообще говоря, лишь кусочно-гладкими системами. Имеющиеся к настоящему времени результаты В.~Ф.~Лазуткина \cite{Lazutkin}, в частности, не дают ответа на вопрос о существовании гладкой и симплектической структур в прообразе окрестности невыпуклых ребер склейки комплекса.
	
	Во-вторых, остается открытым следующий вопрос о топологии таких слоений. Рассмотрим биллиардную книжку $\Omega$, все листы которой --- элементарные плоские биллиарды с выпуклыми углами границы стола (или граница углов не имеет). Пусть перестановки на ребрах (1-клетках комплекса), инцидентных одной и той же вершине, коммутируют. Является ли это достаточным условием лиувиллевой эквивалентности биллиарда на такой книжке $\Omega$ и некоторой гладкой невырожденной ИГС? Т.е. верно ли, что регулярные слои являются торами, а особенности имеют вид боттовских 3-атомов?
	
	Этот вопрос нетривиален: В.~Драговичем и М.~Раднович \cite{Pseudo} изучался класс софокусных биллиардов, названных псевдо-интегрируемыми биллиардами. Вдоль их траекторий сохраняется значение функции $\Lambda$, но регулярные поверхности уровня в $Q^3_h$ гомеоморфны сферам с ручками и проколами. Это обусловлено наличием невыпуклых углов $3\pi/2$ в границе области. На языке биллиардных книжек это означает некоммутативность перестановок $(1)(2\,3)$ и $(1\,2)(3)$ в окрестности общей вершины трех элементарных биллиардов $1, 2, 3$ с выпуклыми углами. Для плоских биллиардов c невыпуклыми углами структура особого слоя и аналог понятия 3-атома изучался В.~А.~Москвиным \cite{Moskvin1}.
	
	В настоящей работе получено два следующих результата. Во-первых, для любого боттовского 3-атома слоения Лиувилля произвольной биллиардной книжки в пункте \ref{s:4:1:flow} показано (утверждение \ref{Ass:1}), что движение частицы в отсутствие потенциала задает \textit{одинаковое} направление на его особых окружностях. Известно \cite{BolFom}, что интегрируемые системы, у которых такое свойство нарушается, не являются топологически устойчивыми (и при построении инварианта Фоменко-Цишанга не рассматривались).
	
	Во-вторых, в пункте \ref{s:4:2:split} показано, что даже достаточно простая биллиардная система с потенциалом может обладать расщепляющимися особенностями (для которых не выполнено \textit{условие нерасщепляемости} Н.~Т.~Зунга, подробнее см. \cite{Zung1, BolOsh}). Эти особенности не содержат точек ранга 0, и все их особые точки --- боттовские. Вместе с тем, слоение Лиувилля в их 4-окрестности нетривиально: в проообразе точки пересечения двух дуг $\Sigma$ имеется особый слой, отличающийся от особых слов атомов на дугах диаграммы.
	
	Как показывает опыт изучения систем физики и классической механики, расщепляющиеся особенности весьма редко встречаются в них. Тем не менее, расщепляющиеся особенности были обнаружены в геодезических потоках на 2-поверхностях вращения в потенциальном поле, изученных Е.~А.~Кантонистовой в \cite{Kantonistova}. Изучение систем с магнитным полем на 2-поверхностях вращения было проведено Е.~А.~Кудрявцевой и А.~А.~Ошемковым в \cite{KudrOsh}. Также расщепляющиеся особенности встречаются в системах с периодическим линейным интегралом на алгебре Ли $e(3)$, изучавшихся в работе И.~К.~Козлова и А.~А.~Ошемкова \cite{KozOshLinear}.
	

\subsection{Задание 3-атома с помощью f-графа}\label{s:1:atoms}

	Переформулировка 2-атомов в виде $f$-графов принадлежит А.\,А.~Ошемкову (см.~\cite{Oshemkov}). Ниже приводится естественное обобщение этой формулировки на случай ориентированных 3-атомов. В случае 3-атомов без звездочек формулировка совпадает с формулировкой для 2-атомов. В случае 3-атомов со звездочками на $ f $-графе необходимо ввести дополнительную структуру, что мы и сделаем ниже.
	
	Заметим также, что в книге А.\,В.~Болсинова, А.\,Т.~Фоменко~\cite{BolFom} уже было введено понятие $ f $-графа со звездочками, являющееся обобщением $ f $-графов на случай 3-атомов. Понятие $ f $-графов со звездочками и понятие $ f $-графов для 3-атомов, изложенное ниже, связаны между собой и друг из друга выводятся. Но в данной работе удобнее использовать понятие $ f $-графов для 3-атомов, введенное нами ниже.

	\begin{example}[построения $f$-графа по 3-атому~$ C_2 $] \label{ex_f-graph_C2}
		Покажем на примере 3-атома~$ C_2 $, изображенного на рис.~\ref{fig_f_graph}, как построить $ f $-граф. Поскольку 3-атом~$ C_2 $ не имеет звездочек, то он представляется в виде прямого произведения 2-атома~$ C_2 $ на окружность. Рассмотрим гладкое погружение 2-атома~$ C_2 $ в плоскость~$ \mathbb{R}^2 $. Такое погружение существует для любого ориентируемого 2-атома согласно теореме А.\,Т.~Фоменко (см.~\cite{BolFom}). Введем граф~$ FC_2 $, который назовем $ f $-графом, соответствующим 3-атому $ C_2 $, следующим образом. Рассмотрим отрезки на 2-атоме~$ C_2 $, соединяющие границы отрицательных колец 2-атома и его критические точки (см. рис.~\ref{fig_f_graph}). Каждая пара таких отрезков, входящих в критическую точку, образует неориентированное ребро $ f $-графа~$ FC_2 $. Вершинами графа~$ FC_2 $ будут концы этих отрезков, лежащих на границе отрицательных колец. Фиксируем ориентацию на критическом уровне 2-атома~$ C_2 $ так, чтобы отрицательные кольца оставались слева. Распространим эту ориентацию на ближайшие некритические уровни. Заметим, что если ориентация задана на одном кольце, то поскольку 2-атом связен, ориентация на всем 2-атоме возникает автоматически. Таким образом, на каждой окружности положительных и отрицательных колец задана ориентация. Эта ориентация на окружностях отрицательных колец дает ориентированные ребра на графе $ FC_2 $. Получившийся граф является $ f $-графом, соответствующим 3-атому~$ C_2 $ без звездочек. В статье~\cite{Oshemkov} А.\,А.~Ошемкова (см. главу 2) показано, что $ f $-графы полностью определяют 2-атом, что эквивалентно утверждению о том, что $ f $-графы полностью определяют 3-атом без звездочек.
		\begin{figure}
			\includegraphics[width=\linewidth]{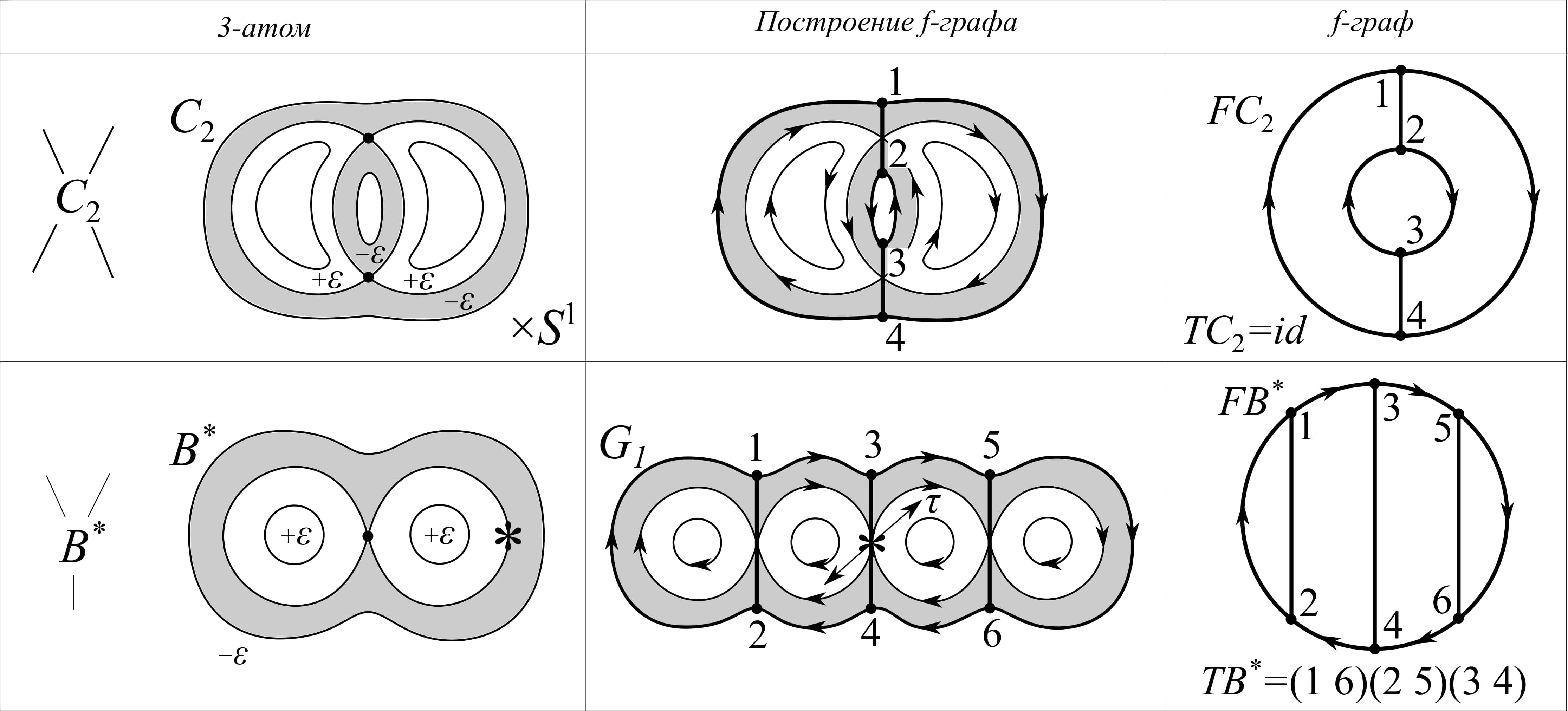}
			\caption{Пример построения $f$-графов~$ FC_2 $, $ FB^* $ по 3-атомам~$ C_2 $, $ B^* $ соответственно.}
			\label{fig_f_graph}
		\end{figure}
	\end{example}
	
	В случае, когда дан 3-атом со звездочками, задания только графа недостаточно для того, чтобы определить 3-атом. Поэтому воспользуемся дублями 3-атомов (см. опр. дубля в~\cite{BolFom}). Дубль с заданной на нем инволюцией полностью определяет 3-атом со звездочкой. В частности, можем считать, что если инволюция тождественна, то задан 3-атом без звездочек.
	\begin{example}[построения $f$-графа по 3-атому~$ B^* $] \label{ex_f-graph_Bz}
		Покажем на примере 3-атома~$ B^* $, как построить $ f $-граф, и дополним его необходимой структурой, чтобы 3-атом со звездочкой однозначно восстанавливался. Рассмотрим дубль 3-атома $ B^* $ --- 2-атом~$ G_1 $ с инволюцией~$ \tau $, изображенной на рис.~\ref{fig_f_graph} как центральная симметрия. Построим аналогично предыдущему примеру $ f $-граф~$ FB^* $ для 2-атома~$ G_1 $ (см. рис.~\ref{fig_f_graph}). Инволюция~$ \tau $ на 2-атоме~$ G_1 $ порождает инволюцию~$ TB^* $ на $ f $-графе~$ FB^* $, являющуюся его автоморфизмом. Инволюция-автоморфизм~$ TB^* $ переводит первую вершину графа в шестую, вторую в пятую, третью в четвертую и наоборот, то есть задается перестановкой $ (1 \, 6) (2 \, 5) (3 \, 4)$. Таким образом, граф~$ FB^* $ вместе с инволюцией~$ TB^* $, также как и дубль~$ G_1 $ вместе с инволюцией~$ \tau $, однозначно задает 3-атом~$ B^* $. Для простоты будем считать, что на $ f $-графе всегда задана инволюция-автоморфизм~$ T $. Если она тождественна, то $ f $-граф задает 3-атом без звездочек. В противном случае $ f $-граф вместе с инволюцией-автоморфизмом~$ T $ задает дубль для 3-атома со звездочками с инволюцией на нем, по которым строится уже сам 3-атом.
	\end{example}
	
	\begin{remark}
		Обратим внимание, что для определения $ f $-атомов и для дальнейшей работы важно, какие кольца у фиксированного 2-атома (3-атома) верхние (лежат на уровне выше критического), а какие нижние (лежат на уровне ниже критического). Иными словами, для каждого из 2-атомов (3-атомов) можно рассмотреть замену функции Морса $ f $ на $ -f $. При этом 2-атом (3-атом) ``переворачивается'' и получается тот же атом, но с другой функцией, то есть перестраивающий окружности (торы) в обратном порядке. Нам важно различать такие атомы. См. пример атома $ B $ и его ``переворота'' на рис.~\ref{fig_f_atoms}. В этом случае получаются разные $ f $-графы.
		\begin{figure}[!ht]
			\centering
			\includegraphics[width=\linewidth]{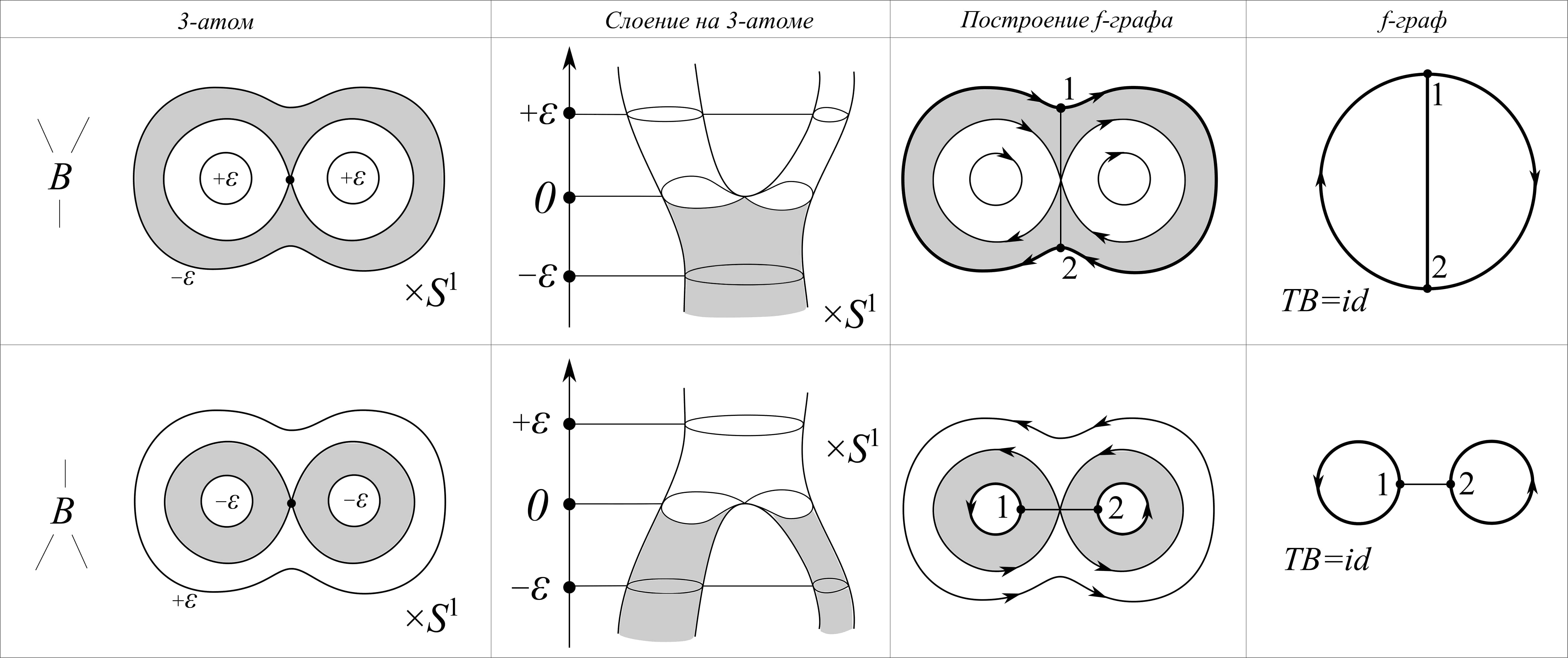}
			\caption[]{Пример двух вариантов перестроек 2-атома~$ B $ и соответствующих им $ f $-графов. На одной строке таблицы 2-атом~$ B $ перестраивает один тор в два. На другой --- 2-атом~$ B $, перестраивающий два тора в один.}
			\label{fig_f_atoms}
		\end{figure}
	\end{remark}
	
	Дадим теперь общее определение $ f $-графа, отвечающее произвольному 3-атому (не только 3-атомам $ C_2 $, $ B^* $ и $ B $).
	\begin{definition}
		Конечный связный граф $ F $ назовем \textit{$ f $-графом}, если он удовлетворяет следующим условиям:
		\begin{enumerate}
			\item
			К каждой вершине графа~$ F $ примыкает ровно одно неориентированное ребро и два ориентированных, из которых одно входит в вершину, а другое выходит из нее. Причем эта вершина может быть началом и концом одного и того же ориентированного ребра в случае, если ориентированное ребро является петлей.
			\item
			Задан автоморфизм $ T $ графа~$ F $, являющийся инволюцией.
		\end{enumerate}
	\end{definition}
	
	Как было упомянуто выше, по любому $ f $-графу однозначно восстанавливается некоторый 3-атом. Однако заметим, что для одного 3-атома может существовать несколько $ f $-графов, которые ему соответствуют. Эта неоднозначность возникает тогда, когда появляются звездочки на 3-атомах, потому что дубли для них определены неоднозначно. Общий алгоритм построения $ f $-графа по произвольному атому и построения атома по произвольному $ f $-графу см. в статье~\cite{Oshemkov} А.\,А.~Ошемкова.
	Также заметим, что неориентированные ребра $ f $-графа --- это сепаратрисы градиентного потока функции Морса на атоме, то есть траектории, входящие или исходящие из критической точки.
	
	Теперь вместо 3-атомов будем оперировать их эквивалентами, $ f $-графами. Оказывается, это удобно для реализации гамильтоновых систем с помощью биллиардов.
	
	\section{Наглядная реализация интегрируемых гамильтоновых систем биллиардными книжками}
	\subsection{Реализация f-графов биллиардными книжками}\label{s:2:1:graphs}	
	\begin{theorem}[В.~В.~Ведюшкина - И.~С.~Харчева~\cite{VedKha}]  \label{th1}
		Гипотеза~\ref{global_hypothesis} Фоменко~$ \mathbf{A} $ верна, а именно, для любого 3-атома (со звездочками или без) алгоритмически строится биллиардная книжка, которая при помощи $ f $-графов реализует данный 3-атом в следующем смысле. В изоэнергетической поверхности~$ Q^3 $ этой биллиардной книжки возникает слоение Лиувилля. Оказывается, что расслоенный прообраз окрестности особого значения интеграла~$ \Lambda $, отвечающего траекториям биллиарда, направленным к или от одного из фокусов, послойно гомеоморфен данному 3-атому.
	\end{theorem}
	
	В статье~\cite{VedKha} предъявлен явный алгоритм построения соответствующей биллиардной книжки на языке 3-атомов. Изложим в данном параграфе этот алгоритм на языке $ f $-графов. Еще раз отметим, что этот алгоритм на языке $ f $-графа 3-атома оказывается существенно проще, чем конструирование биллиарда исходя из самого 3-атома. Поэтому целью данного параграфа является наглядное изложение этого алгоритма. Доказательство того факта, что биллиардная книжка, построенная по этому алгоритму, действительно реализует данный 3-атом, см. в статье~\cite{VedKha}. Итак, проиллюстрируем этот алгоритм на примере $ f $-графа, соответствующего 3-атому $ B^* $.
	
	\begin{figure}
		\centering{\includegraphics[width=0.45\linewidth]{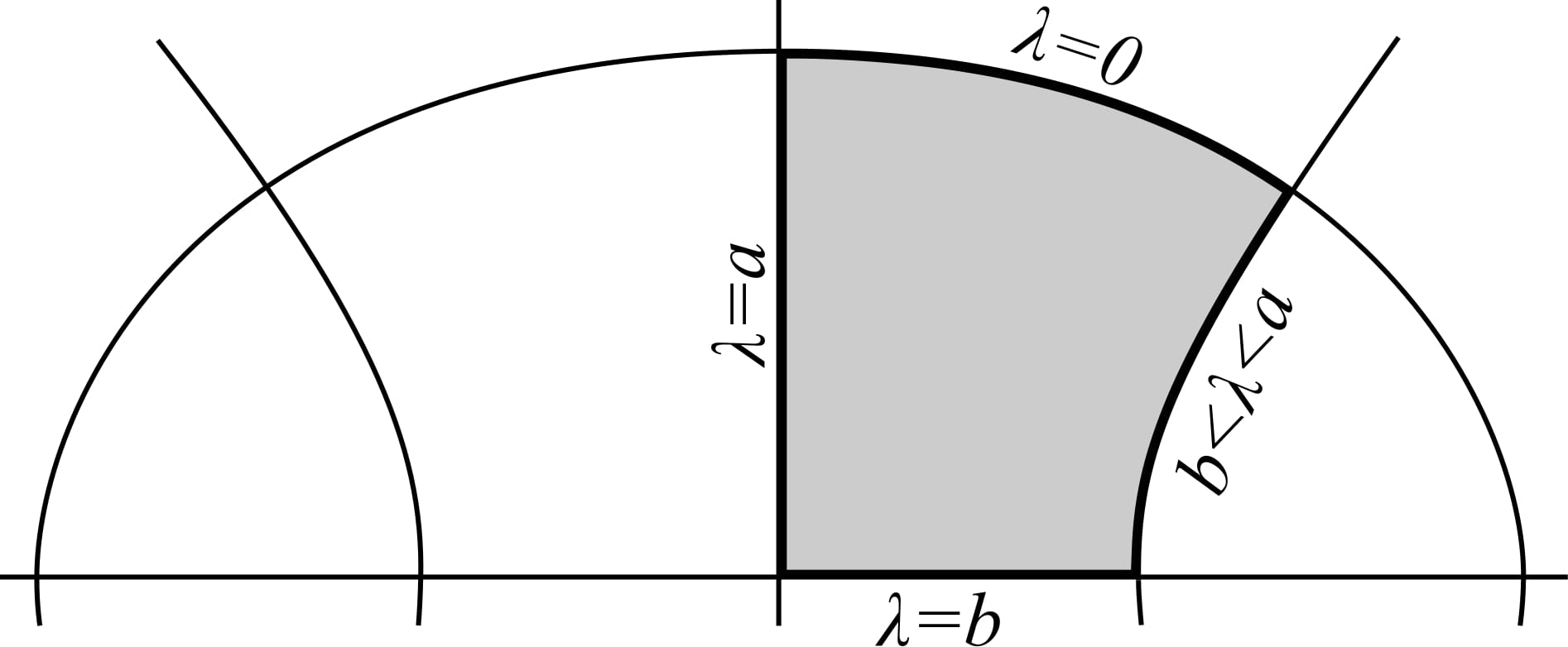}}
		\caption{Область в плоскости, ограниченная дугами софокусных квадрик и являющаяся листами биллиардной книжки, реализующей $ f $-граф. На рисунке обозначены значения, которые принимает параметр $ \lambda $ для каждой квадрики из софокусного семейства~(\ref{eq_confocal}), где $ a $ и $ b $ --- параметры этого семейства.}
		\label{fig_page}
	\end{figure}
	
	\begin{example}[реализации $ f $-графа, соответствующего 3-атому $ B^* $, при помощи биллиардных книжек]
		Напомним, что $ f $-граф~$ FB^* $, соответствующий 3-атому $ B^* $, был построен в примере~\ref{ex_f-graph_Bz} и изображен на рис.~\ref{fig_Bz_book}. Каждой вершине графа~$ FB^* $ сопоставим лист биллиардной книжки, изображенный на рис.~\ref{fig_page}. Он ограничен сверху эллипсом, снизу горизонтальной фокальной прямой, слева вертикальной прямой, проходящей через центр эллипса, а справа гиперболой из софокусного семейства квадрик. Согласно алгоритму, в биллиардной книжке~$ \mathscr{B}B^* $, моделирующей $ f $-граф~$ FB^* $, будет шесть листов. Для упрощения изложения занумеруем произвольным образом вершины и тем самым зададим нумерацию на листах биллиардной книжки.
		
		Выпишем четыре перестановки, на основе которых мы склеим листы биллиардной книжки~$ \mathscr{B}B^* $ друг с другом по правой, левой, нижней и верхней границам.
		Правая перестановка на гиперболе всегда тождественная. Левая перестановка на вертикальной прямой равна инволюции-автоморфизму на $ f $-графе, то есть в случае 3-атома $ B^* $ она равна $ (1 \, 6) (2 \, 5) (3 \, 4) $. В случае любого 3-атома без звездочек, она была бы равна тождественной перестановке. Нижняя перестановка на фокальной прямой является произведением транспозиций, которые задаются на основе неориентированных ребер. В $ f $-графе~$ FB^* $ содержится три неориентированных ребра $ \{1, 2\} $, $ \{3, 4\} $, $ \{5, 6\} $ (запись неориентированного ребра $ \{\alpha, \beta\} $ означает, что это ребро соединяет вершины $ \alpha $ и $ \beta $). Согласно алгоритму, в этом примере нижняя перестановка будет композицией трех транспозиций: $ (1 \, 2) (3 \, 4) (5 \, 6)$.
		Верхняя перестановка на эллипсе пишется на основе простых циклов, образуемых ориентированными ребрами $ f $-графа. Напомним, что в теории графов простой цикл --- это замкнутый обход без повторного прохода по ребру и посещения вершины дважды, за исключением начальной и конечной вершин. Проходим по указанным выше циклам на графе и последовательно выписываем номера вершин, встречающихся вдоль этих циклов, в циклическую перестановку до тех пор, пока не сделаем полный оборот. Один такой цикл дает одну циклическую перестановку. Затем берем композицию всех таких циклических перестановок в произвольном порядке.
		В графе~$ FB^* $ есть один такой цикл. Он образуется ориентированными ребрами $ (1, 3) $, $ (3, 5) $, $ (5, 6) $, $ (6, 4) $, $ (4, 2) $ и $ (2, 1) $ (запись ориентированного ребра $ (\alpha, \beta) $ означает, что ребро начинается в вершине $ \alpha $ и ведет в вершину $ \beta $). В этом случае верхняя перестановка равна $ (1 \, 3 \, 5 \, 6 \, 4 \, 2) $.
		
		Таким образом, биллиардная книжка $ \mathscr{B}B^* $, моделирующая $ f $-граф~$ FB^* $, состоит из шести листов, которые склеены снизу, сверху и слева друг с другом следующим образом. Снизу склеен первый лист со вторым, третий с четвертым, пятый с шестым. Три корешка книжки, возникшие при этой склейке, снабжены транспозициями $ (1 \, 2) $, $ (3 \, 4) $ и $ (5 \, 6) $. Корешки с перестановками $ (1 \, 2) $, $ (3 \, 4) $, $ (5 \, 6) $ соответствуют неориентированным ребрам $ \{1, 2\} $, $ \{3, 4\} $, $ \{5, 6\} $ $ f $-графа~$ FB^* $. Сверху все шесть листов склеены между собой. Образовавшийся корешок снабжен циклической перестановкой $ (1 \, 3 \, 5 \, 6 \, 4 \, 2) $. Слева листы склеены на основе транспозиции $ (1 \, 6) (2 \, 5) (3 \, 4) $. То есть первый лист склеен с шестым, второй с пятым, третий с четвертым. Образовавшиеся три корешка снабжены перестановками $ (1 \, 6) $, $ (2 \, 5) $ и $ (3 \, 4) $. Получившаяся биллиардная книжка~$ \mathscr{B}B^* $ изображена на рис.~\ref{fig_Bz_book}. Оказывается, из доказательства теоремы~\ref{th1} (см.~\cite{VedKha}) вытекает, что окрестность критического уровня интеграла $ \Lambda = b $ в этой книжке послойно гомеоморфна 3-атому~$ B^* $.
		\begin{figure}
			\centering{\includegraphics[width=0.85\linewidth]{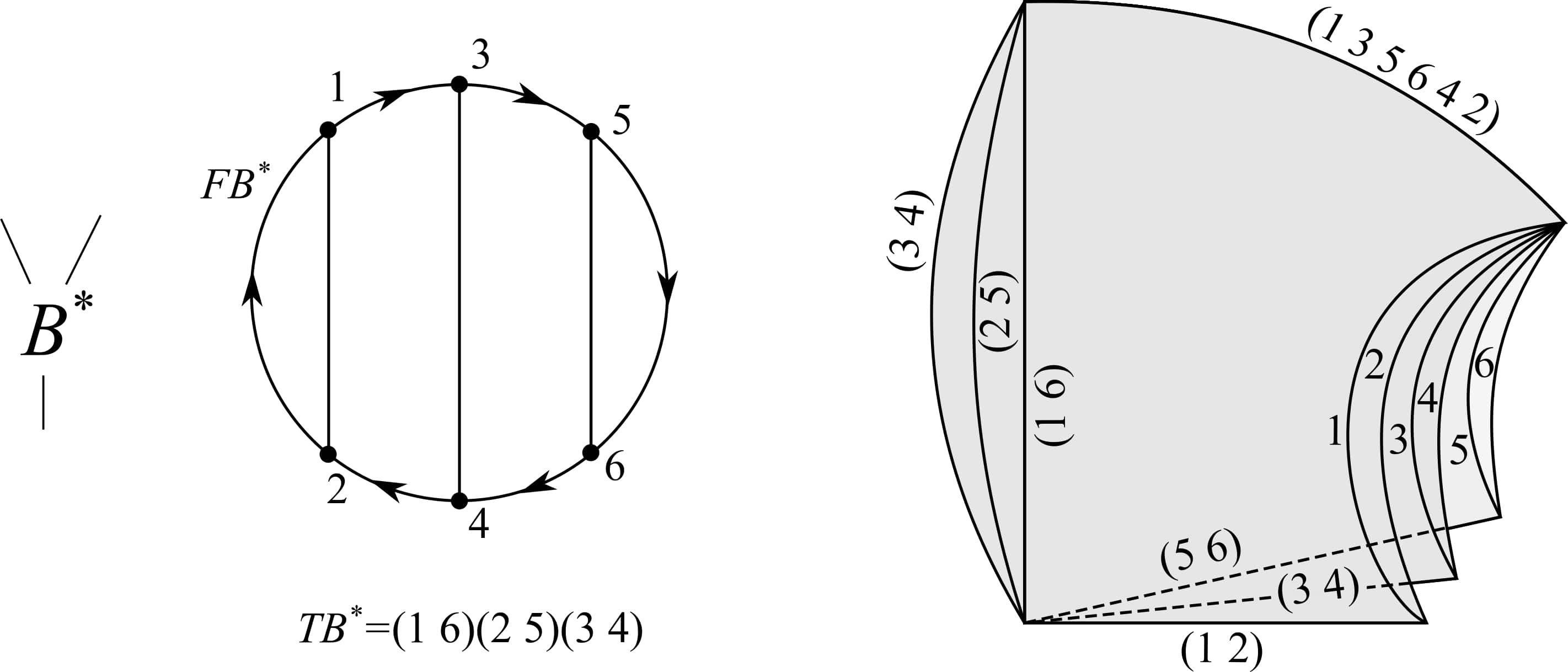}}
			\caption{Пример биллиардной книжки~$ \mathscr{B}B^* $, реализующей $ f $-граф $ FB^* $, иллюстрирующий теорему~\ref{th1}.}
			\label{fig_Bz_book}
		\end{figure}
	\end{example}
	
	В общем случае алгоритм выглядит следующим образом.
	
	\begin{algorithm}[реализации $ f $-графов биллиардными книжками] \label{alg_fgraph}
		Пусть задан произвольный $ f $-граф~$ F $. Этот $ f $-граф соответствует некоторому 3-атому. В результате алгоритма мы собираемся построить такую биллиардную книжку~$ \mathscr{B} $, чтобы в ней встречался этот 3-атом.
		
		Каждой вершине графа~$ F $ сопоставим лист биллиардной книжки, изображенный на рис.~\ref{fig_page}. Он ограничен сверху эллипсом, снизу горизонтальной фокальной прямой, слева вертикальной прямой, проходящей через центр эллипсов, а справа гиперболой из софокусного семейства квадрик. В биллиардной книжке~$ \mathscr{B} $ будет столько листов, сколько вершин в графе~$ F $.  Для упрощения изложения занумеруем произвольным образом вершины, и тем самым зададим нумерацию на листах биллиардной книжки.
		
		Выпишем четыре перестановки, на основе которых мы склеим листы биллиардной книжки~$ \mathscr{B} $ друг с другом по правой, левой, нижней и верхней границам.
		Правая перестановка на гиперболе всегда тождественная и никакие два листа не склеены друг с другом по правой границе. Левая перестановка на вертикальной прямой равна инволюции-автоморфизму на $ f $-графе. Заметим, что эта перестановка тождественная тогда и только тогда, когда $ f $-граф соответствует 3-атому без звездочек.
		Нижняя перестановка на фокальной прямой задается на основе неориентированных ребер графа $ F $ следующим образом. Если дано неориентированное ребро $ \{\alpha, \beta\} $, соединяющее вершины $ \alpha $ и $ \beta $, то рассмотрим транспозицию $ (\alpha \, \beta) $. Перемножим все такие транспозиции в произвольном порядке и получим нижнюю перестановку.
		Верхняя перестановка на эллипсе пишется на основе простых циклов, образуемых ориентированными ребрами $ f $-графа. Проходим по каждому из таких циклов на графе и последовательно выписываем номера вершин, встречающихся вдоль этих циклов, в циклическую перестановку до тех пор, пока не сделаем полный оборот. Один такой цикл дает одну циклическую перестановку. Затем берем композицию всех таких циклических перестановок в произвольном порядке.
		
		Далее опишем механизм склейки листов книжек на основе перестановок. Пусть дана верхняя перестановка. Опишем как склеить занумерованные листы по верхней границе. Раскладываем верхнюю перестановку в произведение независимых циклов (она уже будет разложена в такое произведение, если она строилась согласно действиям, описанным выше).
		Рассмотрим один из циклов $ (\alpha_1 \, \alpha_2 \, ... \, \alpha_k) $ некоторой длины $ k $ из этого разложения. Склеим друг с другом листы, принадлежащие орбите этого цикла, то есть листы $ \alpha_1, \alpha_2, ..., \alpha_k $, и только их. При этом образуется корешок. Снабдим этот новый корешок циклической перестановкой $ (\alpha_1 \, \alpha_2 \, ... \, \alpha_k) $. Такие же действия выполним для каждого из циклов и для всех сторон листов: верхней, нижней, левой и правой. В случае 3-атомов без звездочек, такая склейка произойдет только по верхней и нижней границам, в случае 3-атомов со звездочками --- по верхней, нижней и левой границам. Получившаяся биллиардная книжка~$ \mathscr{B} $ реализует заданный изначально $ f $-граф~$ F $.
		
		Еще раз напомним, что в данной статье мы наглядно описываем только алгоритм действий. Доказательство того, что этот алгоритм действительно реализует заданный $ f $-граф дано в работе~\cite{VedKha}.
	\end{algorithm}
	
	\begin{remark}
		Заметим, что при задании нижней и верхней перестановок в алгоритме~\ref{alg_fgraph} композицию можно брать в произвольном порядке. Действительно, нижняя перестановка на фокальной прямой задается на основе неориентированных ребер $ f $-графа (см. алгоритм~\ref{alg_fgraph}). Согласно определению $ f $-графа, из каждой вершины выходит ровно одно неориентированное ребро. Следовательно, каждая из вершин встречается только в одной транспозиции. Значит, числа в транспозициях не пересекаются, транспозиции независимы, и их можно перемножать в произвольном порядке.
		
		Верхняя перестановка на эллипсе пишется на основе простых циклов, образуемых ориентированными ребами $ f $-графа (см. алгоритм~\ref{alg_fgraph}). Из определения $ f $-графа следует, что в каждую вершину графа входит и из каждой вершины выходит ровно одно ориентированное ребро. Поскольку при построении верхней перестановки рассматриваются простые циклы, то каждую из вершин может содержать только один из таких циклов. Значит, каждая из вершин встретится только в одной циклической перестановке, эти циклы независимы, их орбиты не пересекаются, и их можно перемножать в произвольном порядке.
	\end{remark}
	
	\subsection{Пример. Наглядная реализация биллиардами случая Жуковского динамики твердого тела на одном из уровней энергии} \label{s:2:molec}
	В статье~\cite{VedKha2} доказана следующая теорема.
	\begin{theorem}[В.~В.~Ведюшкина-И.~С.~Харчева] \label{th2}
		Гипотеза~\ref{global_hypothesis} Фоменко~$ \mathbf{B} $ верна, то есть любая грубая молекула реализуется биллиардными книжками. Более точно: по любой грубой молекуле алгоритмически строится биллиардная книжка класса~$ \mathbf{b} $ с каноническим квадратичным интегралом~$ \Lambda $, отвечающим параметру каустики, такая что грубая молекула, соответствующая этой системе, изоморфна заданной изначально грубой молекуле.
	\end{theorem}
	
	В данном параграфе приведем наглядную реализацию грубой молекулы с помощью биллиардных книжек для случая Н.~Е.~Жуковского динамики твердого тела. Напомним, как выглядит эта система.
	
	Классические уравнения Эйлера-Пуассона, описывающие движение твердого тела с закрепленной точкой в поле силы тяжести, в системе координат, оси которых направлены вдоль главных осей инерции тела, имеют следующий вид (см., например,~\cite{Kozlov}).
	\begin{equation}\label{eq_Euler_Poiss}
		\begin{gathered}
			A \dot{\omega} = A \omega \times \omega - P r \times \nu, \\
			\dot{\nu} = \nu \times \omega.
		\end{gathered}
	\end{equation}
	Фазовые переменные здесь таковы: $ \omega $ -- вектор угловой скорости, $ \nu $ -- единичный вертикальный вектор. Параметрами системы являются: диагональная матрица $ A = diag(A_1, A_2, A_3) $, задающая тензор инерции твердого тела, $ P $ --- вес тела, $ r $ --- вектор с началом в неподвижной точке тела и концом в центре масс тела. Запись $ a \times b $ означает векторное произведение $ a $ и $ b $ в $ R^3 $.
	
	Вектор $ A\omega $ имеет смысл кинетического момента твердого тела относительно неподвижной точки. Н.~Е.~Жуковский исследовал задачу о движении твердого тела, имеющего полости, целиком заполненные идеальной несжимаемой жидкостью, совершающей безвихревое движение. В этом случае кинетический момент тела равен $ A\omega + \lambda $. Здесь $ \lambda $ --- постоянный (в системе координат, связанной с телом) вектор, характеризующий циклические движения жидкости в полостях. Аналогичный вид кинетический момент тела имеет в случае, когда в теле закреплен маховик, ось которого направлена вдоль вектора $ \lambda $. Такую механическую систему называют гиростатом. Движение гиростата в поле силы тяжести описываются системой уравнений.
	\begin{equation}\label{eq_Zukovsky}
		\begin{gathered}
			A \dot{\omega} = (A \omega + \lambda) \times \omega - P r \times \nu, \\
			\dot{\nu} = \nu \times \omega,
		\end{gathered}
	\end{equation}
	частным случаем которой при $ \lambda = 0 $ является система~\ref{eq_Euler_Poiss}. У системы~\ref{eq_Zukovsky} всегда существуют геометрический интеграл
	$$ F = \langle\nu, \nu\rangle = 1,$$
	интеграл энергии
	$$ E = \frac{1}{2} \langle A\omega, \omega \rangle + P\langle r, \nu \rangle,$$
	и интеграл площадей
	$$ G = \langle A\omega + \lambda, \nu \rangle.$$
	Запись $ \langle a, b \rangle $ означает скалярное произведение произведение $ a $ и $ b $ в $ R^3 $.
	
	\begin{figure}
		\centering{\includegraphics[width=0.4\linewidth]{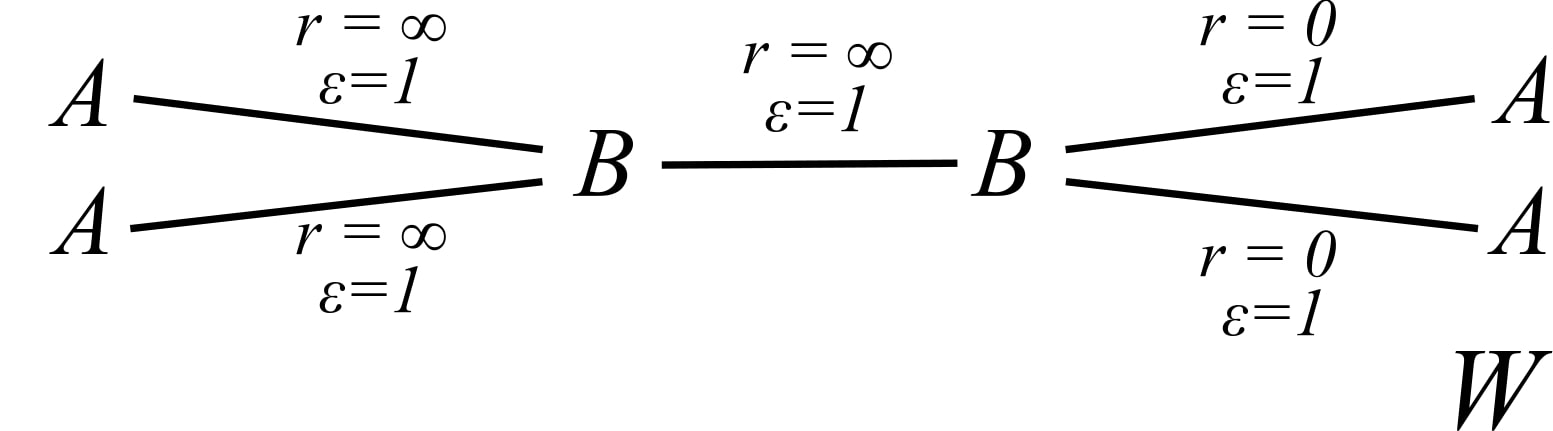}}
		\caption{Меченая молекула~$ W $, соответствующая случаю Жуковского динамики твердого тела на одном из уровней энергии.}
		\label{fig_mol_zukovsky}
	\end{figure}
	
	Оказывается, уравнение~\ref{eq_Zukovsky} является гамильтоновым на совместных четырехмерных поверхностях уровня геометрического интеграла и интеграла площадей (см. доказательсто в книге~\cite{BolFom} т. 2). Для такой системы А.~А.~Ошемков и П.~Й.~Топалов провели лиувиллев анализ (см.~\cite{BolFom} т. 2) и, в частности, вычислили топологические инварианты Фоменко-Цишанга. На одном из уровней энергии возникает молекула~$ W $, изображенная на рис.~\ref{fig_mol_zukovsky}. Используем алгоритм, описанный в работе~\cite{VedKha2} и доказывающий теорему~\ref{th2}, и реализуем эту молекулу с помощью биллиардных книжек.
	
	\begin{figure}
		\centering{\includegraphics[width=0.9\linewidth]{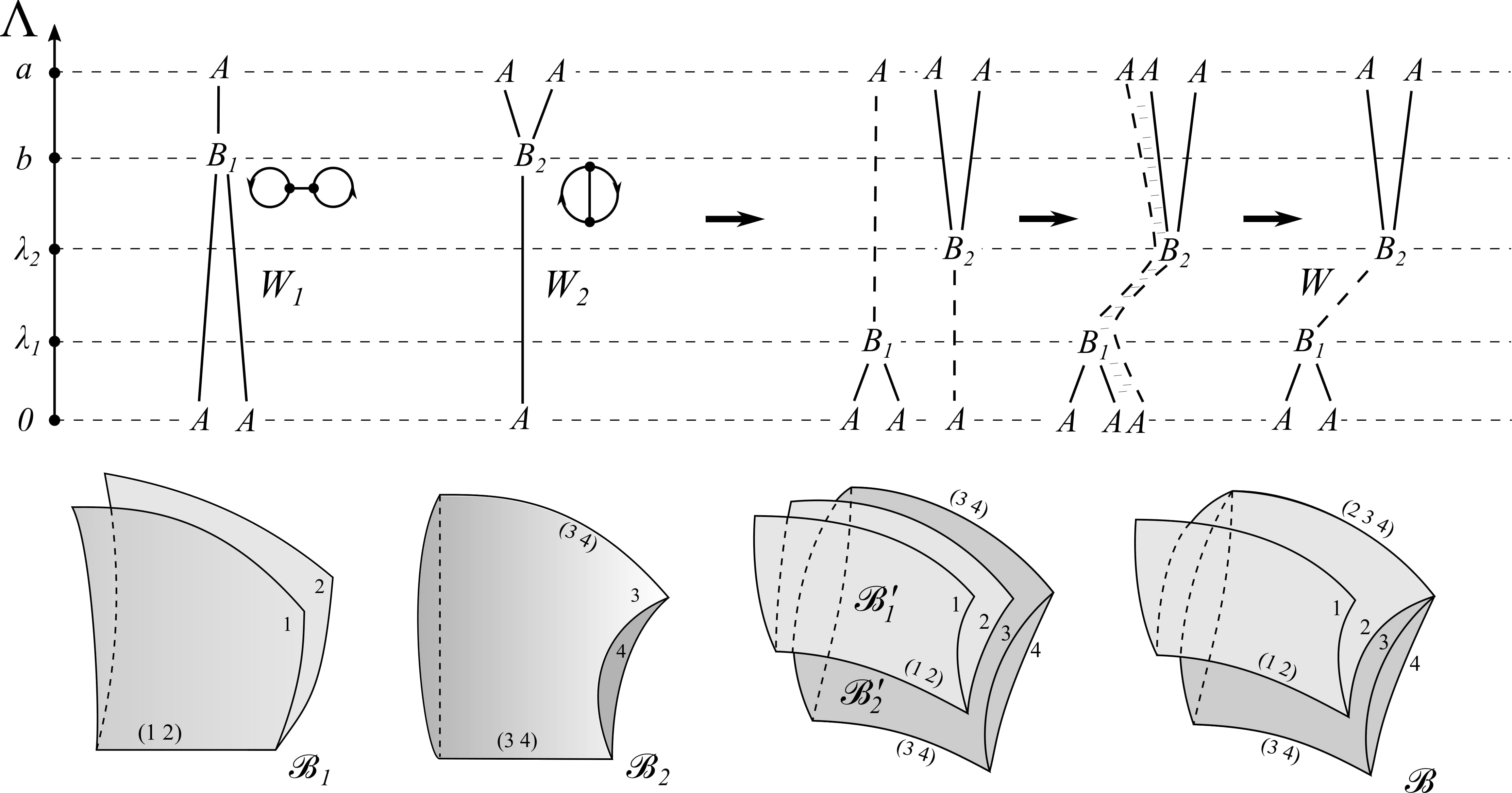}}
		\caption{Построение биллиардной книжки~$ \mathscr{B} $, реализующей грубую молекулу~$ W $.}
		\label{fig_bb_zukovsky}
	\end{figure}
	
	Рассматриваемая молекула~$ W $ содержит два 3-атома $ B $. Один перестраивает два тора в один (для определенности обозначим его через $ B_1 $), другой ($ B_2 $) один тор в два. $ f $-графы для обоих случаев приведены на рис.~\ref{fig_f_atoms}. Используя алгоритм~\ref{alg_fgraph}, построим две биллиардные книжки $ \mathscr{B}_1 $, $ \mathscr{B}_2 $, реализующие эти два $ f $-графа. Эти биллиардные книжки $ \mathscr{B}_1 $, $ \mathscr{B}_2 $, также как и их грубые молекулы $ W_1 $, $ W_2 $, представлены на рис.~\ref{fig_bb_zukovsky}. Напомним, что интеграл $ \Lambda $ является параметром каустики (см.~(\ref{eq_integral})), а параметры $ a $ и $ b $ являются параметрами семейства софокусных квадрик, задаваемых уравнением~(\ref{eq_confocal}). Заметим, что у обеих молекул $ W_1 $ и $ W_2 $ 3-атомы $ B $ находятся на уровне интеграла $ \Lambda = b $. В требуемой молекуле 3-атомы $ B $ находятся на разных уровнях интеграла. Зафиксируем произвольные константы $ \lambda_1 $ и $ \lambda_2 $, такие что $ 0 < \lambda_1 < \lambda_2 < b $, и расположим 3-атом $ B_1 $, перестраивающий два тора в один, на уровне $ \Lambda = \lambda_1 $, 3-атом $ B_2 $, перестраивающий один тор в два, --- на уровне $ \Lambda = \lambda_2 $. Для этого в биллиардных книжках $ \mathscr{B}_1 $ и $ \mathscr{B}_2 $ нужно поменять нижнюю границу. Ранее у биллиардной книжки $ \mathscr{B}_1 $ нижняя граница являлась горизонтальной прямой (квадрикой из семейства с параметром $ b $). Теперь в новой книжке $ \mathscr{B}'_1 $ нижняя граница станет невыпуклой дугой эллипса с параметром $ \lambda_1 $. Аналогично, в новой книжке $ \mathscr{B}'_2 $ нижняя граница станет невыпуклой дугой эллипса с параметром $ \lambda_2 $. Для того, чтобы склеить 3-атомы $ B $ по пунктирным ребрам графа, как показано на рис.~\ref{fig_bb_zukovsky}, склеим биллиардные книжки $ \mathscr{B}'_1 $ и $ \mathscr{B}'_2 $ друг с другом по верхней границе. Для этого приклеим второй лист книжки $ \mathscr{B}'_1 $ к корешку, соединяющему третий и четвертый лист книжки $ \mathscr{B}'_2 $. Полученная биллиардная книжка~$ \mathscr{B} $ представлена на рис.~\ref{fig_bb_zukovsky} и реализует требуемую молекулу~$ W $. Иными словами, грубая молекула, вычисленная для биллиардной книжки~$ \mathscr{B} $, будет иметь вид молекулы~$ W $ без меток.
	
	\begin{theorem}
		Построенная нами биллиардная книжка~$ \mathscr{B} $ (см. рис.~\ref{fig_bb_zukovsky}, биллиард справа), рассматриваемая как интегрируемая система на трехмерном уровне постоянной энергии, грубо лиувиллево эквивалентна системе Жуковского в динамике твердого тела на одном из ее уровней энергии.
	\end{theorem}
	
	\begin{remark}
		Рассмотрим подробнее, что произошло, когда мы склеили две биллиардные книжки $ \mathscr{B}'_1 $ и $ \mathscr{B}'_2 $ между собой: как поменялось движение материальной точки на биллиардной книжке на каждом уровне интеграла $ \Lambda $, как это изменение отразилось на слоях и соответственно на грубой молекуле (см. рис.~\ref{fig_bb_movement}). Для того, чтобы описывать такое движение, введем обозначения. Запись $ \uparrow_{\alpha}$ (запись $ \downarrow_{\alpha}$) означает, что материальная точка идет по листу с номером $ \alpha $ наверх (соответственно вниз). Последовательная запись $ \uparrow_{\alpha} \downarrow_{\beta} $ означает, что материальная точка идет сначала по листу $ \alpha $ наверх, отражается от верхней границы, затем по листу $ \beta $ вниз, отражается от нижней границы, и после этого движение замыкается в том смысле, что оно снова идет по листу $ \alpha $ наверх. Такая последовательность может состоять из произвольного числа стрелок, не только из двух. Это обозначение с точностью до циклического сдвига однозначно определяет движение по биллиардной книжке.
		\begin{figure}
			\includegraphics[width=\linewidth]{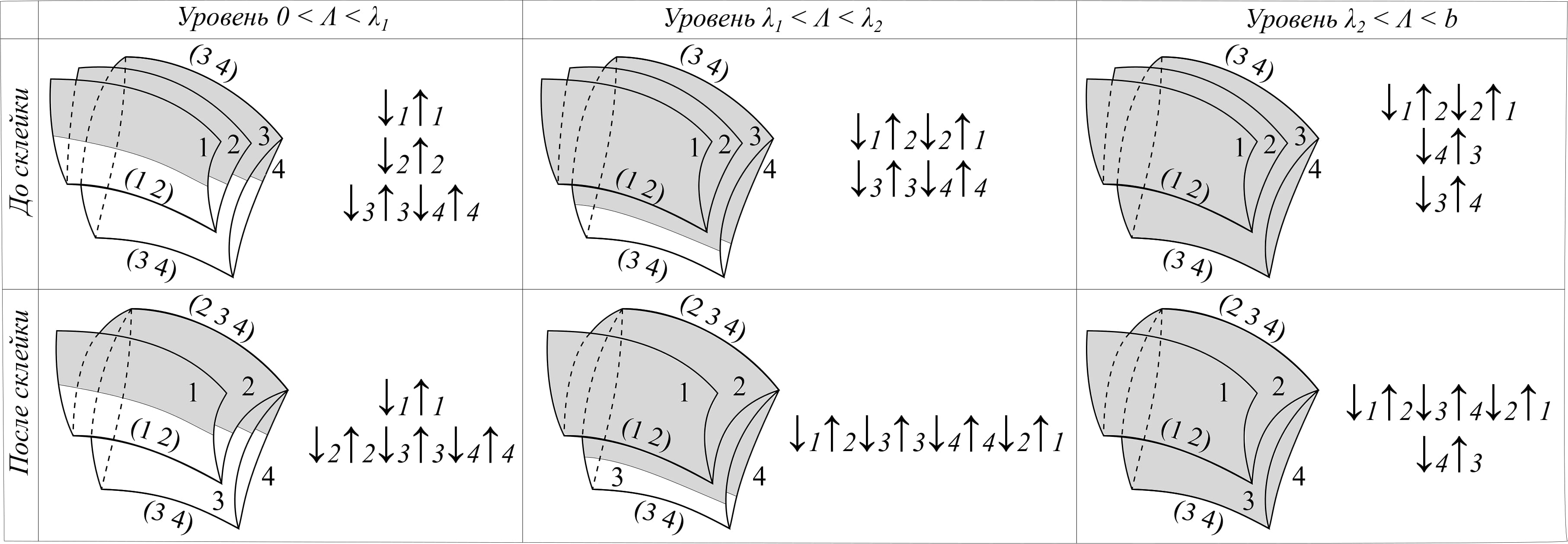}
			\caption{Движение материальной точки на биллиардных книжках на каждом из уровней интеграла $ \Lambda $. На рисунке закрашены области возможного движения.}
			\label{fig_bb_movement}
		\end{figure}
		
		\begin{itemize}
			\item
			Рассмотрим уровень интеграла $ \Lambda = \lambda $, где $ 0 < \lambda < \lambda_1 $. На этом уровне траектория материальной точки всегда направлена по касательной к эллипсу с параметром $ \lambda $. Значит, материальная точка может находиться только выше этого эллипса. Это свойство траекторий на уровне интеграла $ \Lambda = \lambda $ задает область возможного движения, изображенную на рис.~\ref{fig_bb_movement}. До склейки было три непересекающихся между собой движения материальной точки по биллиардной книжки, каждое из которых соответствует некоторому ребру молекулы (тору): $ \downarrow_{1} \uparrow_{1} $, $ \downarrow_{2} \uparrow_{2} $, $ \downarrow_{3} \uparrow_{3} \downarrow_{4} \uparrow_{4} $. После склейки движение материальной точки по второму листу (вверх и вниз) объединилось с движением по третьему и четвертому листу, то есть образовалось два движения: $ \downarrow_{1} \uparrow_{1} $, $ \downarrow_{2} \uparrow_{2} \downarrow_{3} \uparrow_{3} \downarrow_{4} \uparrow_{4} $. Это означает, что тор, соответствующий одному из нижних ребер 3-атома $ B_1 $, приклеился к тору, соответствующему нижнему ребру 3-атома $ B_2 $.
			\item
			Рассмотрим теперь уровень $ \Lambda = \lambda $ посередине, то есть $ \lambda_1 < \lambda < \lambda_2 $. До склейки было два непересекающихся между собой движения: по книжке~$ \mathscr{B}'_1 $ ($ \downarrow_{1} \uparrow_{2} \downarrow_{2} \uparrow_{1} $) и по книжке~$ \mathscr{B}'_2 $ ($ \downarrow_{3} \uparrow_{3} \downarrow_{4} \uparrow_{4} $). После склейки образовалось одно движение $ \downarrow_{1} \uparrow_{2} \downarrow_{3} \uparrow_{3} \downarrow_{4} \uparrow_{4} \downarrow_{2} \uparrow_{1} $, проходящее по всем листам обеих книжек. Это означает, что два тора, которые раньше находились на уровне $ \Lambda = \lambda $ склеились в один.
			\item
			На уровне $ \Lambda = \lambda $, где $ \lambda_2 < \lambda < b $ до склейки было три движения: $ \downarrow_{1} \uparrow_{2} \downarrow_{2} \uparrow_{1} $,  $\downarrow_{4} \uparrow_{3} $, $ \downarrow_{3} \uparrow_{4} $. Первое проходило по книжке $ \mathscr{B}'_1 $. Второе и третье движения материальной точки были направлены по и против часовой стрелки на третьем и четвертом листах книжки $ \mathscr{B}'_2 $. Благодаря склейке первое движение объединилось с одним из последних двух. В итоге, в новой книжке на уровне $ \Lambda = \lambda $ образовалось два движения: $ \downarrow_{1} \uparrow_{2} \downarrow_{3} \uparrow_{4} \downarrow_{2} \uparrow_{1}$ и  $\downarrow_{4} \uparrow_{3} $. Это означает, что тор, соответствующий одному из верхних ребер 3-атома $ B_2 $, приклеился к тору, соответствующему верхнему ребру 3-атома $ B_1 $.
		\end{itemize}
		Таким образом, склейка биллиардных книжек $ \mathscr{B}'_1 $ и $ \mathscr{B}'_2 $ дала склейку двух молекул между собой.
	\end{remark}
	
	\begin{remark}
		В работе~\cite{VedKha2} предъявлен алгоритм, который позволяет реализовывать любую грубую молекулу (не только описанную выше конкретную молекулу $ W $). Он описывает, как нужно склеивать между собой биллиардные книжки по верхней границе в случае произвольных 3-атомов ($ f $-графов), встречающихся в молекуле. Кроме того, молекула может состоять из любого количества произвольных 3-атомов ($ f $-графов) и ребер, соединяющих их. Согласно этому алгоритму, нужно в произвольном порядке приклеивать ребра, соединяющие 3-атомы, выполняя правильным образом переклейку верхней границы листов книжки, до тех пор, пока не получим требуемую грубую молекулу.
	\end{remark}
	
	\begin{theorem}
		Построенная нами биллиардная книжка~$ \mathscr{B} $ (см. рис.~\ref{fig_bb_zukovsky}, биллиард справа), рассматриваемая как интегрируемая система на трехмерном уровне постоянной энергии, грубо лиувиллево эквивалентна системе Жуковского в динамике твердого тела на одном из ее уровней энергии.
	\end{theorem}
	
	Как отметила В.\,В.~Ведюшкина, оказывается, построенная нами биллиардная книжка оказалась настолько ``удачной'', что она реализует не только грубую молекулу системы Жуковского, но реализует также и полный инвариант Фоменко-Цишанга данной гамильтоновой системы. То есть реализует так называемую меченую молекулу (грубую молекулу с числовыми метками). А именно, верна следующая теорема.
	
	\begin{theorem} \label{th_bb_zukovsky}
		Меченая молекула интегрируемой системы, задаваемой  биллиардной книжкой на рис.~\ref{fig_bb_zukovsky} (биллиард справа) изоморфна меченой молекуле (инварианту Фоменко-Цишанга) системы Жуковского на одном из ее трехмерных уровней энергии. Другими словами, эта биллиардная  книжка  ``лиувиллево эквивалентна'' системе Жуковского на соответствующем уровне энергии.
	\end{theorem}
	
	\begin{remark}
		Итак, две совершенно разные по своей природе динамические интегрируемые системы оказались лиувиллево эквивалентными, то есть имеют ``одинаковые'' замыкания интегральных траекторий общего положения.
	\end{remark}
	
	Доказательство теоремы~\ref{th_bb_zukovsky} здесь опущено. Оно будет приведено в другой публикации.

	\section{Числовые инварианты интегрируемых биллиардов}\label{s:3:local}
	
	Приведем вначале формулировку локальной гипотезы А.~Т.~Фоменко, сформулированную в \cite{DAN20}:
	
	\begin{hypothesis}[A.~Т.~Фоменко] Каждое из перечисленных ниже подмножеств меченой молекулы ИГС реализуется как подмножество меченой молекулы некоторого подходящего интегрируемого биллиарда:
		
		\begin{enumerate}
			\item (реберный инвариант) ребро с произвольными выбранными метками $(r, \varepsilon)$;
			\item (усиление пункта 1) ребро с произвольными выбранными метами $(r, \varepsilon)$ между двумя выбранными атомами;
			\item (инвариант семьи) произвольная выбранная метка $n$ на некоторой семье;
			\item (усиление пункта 3) произвольная выбранная метка $n$ на произвольной выбранной семье $S$;
			\item (меченая окрестность семьи) произвольная выбранная семья с произвольной выбранной меткой n на ней и произвольными выбранными реберными инвариантами $(r, \varepsilon)$ на ее внешних ребрах;
			\item (меченая окрестность ребра). Две семьи $S_1$ и $S_2$ с метками $n_1$ и $n_2$, произвольные выбранные граничные торы которых соединены ребром с произвольным выбранным реберным инвариантом $(r, \varepsilon)$.
		\end{enumerate}
	\end{hypothesis}
	
	\begin{remark}
		Пункты 1 и 3 локальной гипотезы посвящены реализации произвольного значения числовым меток $n$ на семье и пары $r, \varepsilon$ на ребре хотя бы одной молекулы. В разделах 2, 4, 5, 6 обсуждается реализация произвольных подграфов Фоменко-Цишанга вблизи ребра или выбранной семьи, т.~е. имеются условия не только на метки, но и на 3-атомы.
		
		Эта задача существенно сложнее. Как ожидается, реализация метки $n$ на произвольной семье, а также совместная реализация наборов реберных инвариантов и метки $n$ потребует развития новых, более тонких методов и конструкций.
	\end{remark}
	
	Напомним, что семья в молекуле --- это связный подграф седловых атомов, на внутренних ребрах которого стоят метки $r = \infty$, а на каждом внешнем ребре, соединяющем атом семьи с эллиптическим атомом $A$, метка $r \ne \infty$. Это значит, что метка $r$ на внешних ребрах семьи принимают произвольное значение из $\mathbb{Q}\,\mathrm{mod}\,1$.
	
	В настоящее время доказана реализация каждого из числовых инвариантов биллиардов независимо от остальных, т.~е. пункты 1, 3 (о значениях меток $n$ и пары $r, \varepsilon$). Этот результат кратко изложен в работе \cite{DAN20} В.~В.~Ведюшкиной, В.~А.~Кибкало и А.~Т.~Фоменко вместе ``атласом'' найденных ранее биллиардов, реализующих определенные продвижения в доказательстве пунктов 2, 4. Реализация метки $n$ (пункта 3) более подробно описана в работе \cite{VedKib20} В.~В.~Ведюшкиной и В.~А.~Кибкало.
	
	\begin{theorem}[см. \cite{DAN20}]
		Для любого значения целочисленной метки $n$ или любой пары значений $r \in \{\mathbb{Q}\, \mathrm{mod}\, 1\} \cup \{\infty\},\, \varepsilon = \pm 1$ построена биллиардная книжка, слоение Лиувилля которой содержит некоторую семью или ребро соответственно с такими числовыми инвариантами.
	\end{theorem}

	\subsection{Реализация реберного инварианта}
	Реализация произвольных пар меток $(r, \varepsilon)$ на некотором ребре выполнена В.~В.~Ведюшкиной (см. \cite{DAN20}) в утверждении \ref{Ass:VedREps} (где полностью доказан пункт 1 гипотезы) и утверждении \ref{Ass:Nomer2}. В нем на основе ряда результатов был составлен ``атлас'' биллиардов, которыми реализована значительная часть возможных случаев (т.~е. значений меток и типов атомов на концах ребра). Как ожидается, ответ в остальных случаях потребует продвижений в общей и нетривиальной проблеме классификации биллиардных книжек.
	
	\begin{assertion}[В.~В.~Ведюшкина]\label{Ass:VedREps}
		Любая пара числовых меток $(r, \varepsilon)$ реализуется как метки на ребре меченой молекулы $W^{*}$ слоения Лиувилля некоторого биллиарда.
	\end{assertion}
	
	\textit{Замечание 1.}
	Ориентация 3-многообразия $Q^3$ влияет на выбор допустимых базисов на граничных торах 3-атомов и на матрицу перехода между ними. Тем самым, метки могут измениться в зависимости от типа 3-атомов на концах ребра (эллиптический или седловой) и конечности метки $r$ (подробнее см. \cite{BolFom}, т.1 раздел 4.5).
	
	Утверждение \ref{Ass:Nomer2} содержит результаты доказательства пункта \textbf{2} (об окрестности некоторого ребра с метками) локальной гипотезы Фоменко:
	
	\begin{assertion}[В.~В.~Ведюшкина]\label{Ass:Nomer2}
		Пункт \textbf{2} локальной гипотезы Фоменко доказан в случаях, указанных в табл.~\ref{tab_c1}. В семи случаях подходящими биллиардами реализован требуемый числовой инвариант $(r,\ \varepsilon)$ на ребре, соединяющем произвольные атомы требуемых типов (эллиптического или седлового). В четырех других случаях  требуемый числовой инвариант реализован на ребре, соединяющем атомы требуемых типов, но седловые атомы не произвольны, а принадлежат конкретным сериям $B_n$ и $C_n$ (в обозначениях \cite{FokDomain}).
	\end{assertion}
	
	\begin{table}[ht]
		\begin{center}
			\begin{tabular}[t]{|c||c|c|c|c|c|c|c|c|c|c|c|c|c|}
				\hline
				метки & $A - A$ & $A - V$ & $V_1 - V_2$  \\
				\hline
				$r = p \slash q,    $ & см. \cite{Lens}, рис. \ref{Fig:REps}a  & $V=B$, рис. \ref{Fig:REps}b  & см. \cite{VedFomGeod} $V_1=C_k,\ V_2=C_n$  \\
				$  \varepsilon = 1$ &    &   &    \\
				\hline
				$r = p\slash q,  $ &см.  \cite{Lens}   & ? &см. \cite{VedFomGeod} $V_1=C_k,\ V_2=C_n$ \\
				$ \varepsilon = -1$ &    &  &  \\
				\hline
				$r = \infty, $ &см. \cite{MagPot}   &см. \cite{VedKha} алгоритм & алгоритм    \\
				$  \varepsilon = 1$ &  и   Предл. 7.1, \cite{VedFomGeod}  &Ведюшкиной-Харчевой  & Ведюшкиной-Харчевой  \\
				&    &  & для грубых молекул \\
				\hline
				$r = \infty, $ &см. \cite{MagPot}  &см. \cite{VedKha}  алгоритм & см.  Предл. 7.2, \cite{VedFomGeod},
				\\
				$ \varepsilon = -1$  &   & Ведюшкиной-Харчевой &  $V_1=V_2=B_n$ \\
				\hline
			\end{tabular}
			\caption{Случаи реализации комбинаций меток  $r$ и $\varepsilon$
				в молекулах Фоменко--Цишанга интегрируемых биллиардов.}
			\label{tab_c1}
		\end{center}
	\end{table}

	В табл.~\ref{tab_c1} для различных вариантов меток и атомов указаны работы \cite{VedKha}, \cite{VedFomGeod}, \cite{Lens}, \cite{MagPot}  (В.~В.~Ведюшкиной, И.~С.~Харчевой и А.~Т.~Фоменко), в которых вводились конструкции биллиардных столов, реализующих указанные сочетания. Отметим также недавний результат С.Е.~Пустовойтова, упомянутый в \cite{MagPot}.
	
\begin{figure}[!htb]
		\centering{\includegraphics[width=0.8\linewidth]{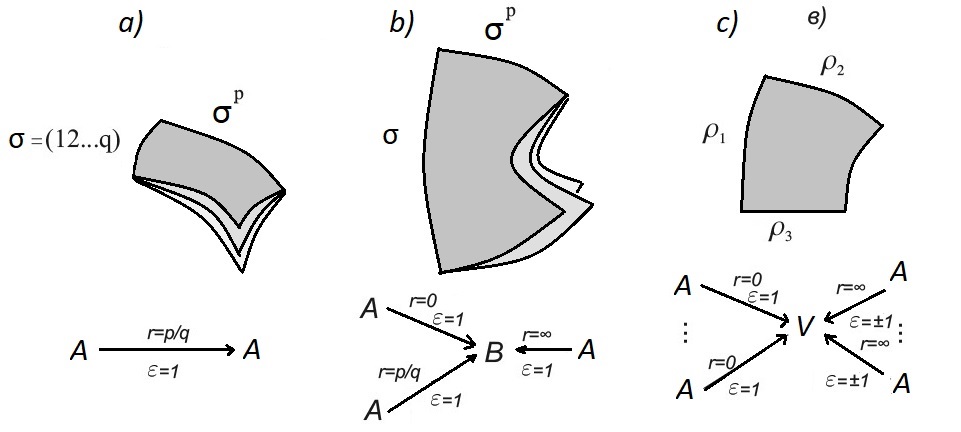}}
		\caption{Биллиарды, реализующие некоторые числовые инварианты из утверждения \ref{Ass:Nomer2}.}
		\label{Fig:REps}
	\end{figure}

	На рис. \ref{Fig:REps} изображены биллиардные столы (и меченые молекулы их слоений Лиувилля), реализующие некоторые случаи утверждения \ref{Ass:Nomer2}: $A \, \cfrac{r = p \slash q}{\varepsilon = \pm 1}\, A$ на рис. \ref{Fig:REps}a, $A\,  \cfrac{r = p \slash q}{\varepsilon = \pm 1}\, A$ на рис. \ref{Fig:REps}b и $A \cfrac{r = \infty}{\varepsilon = \pm 1} V$ с произвольным седловым атомом $V$ на рис. \ref{Fig:REps}c. На рис. \ref{Fig:REps}a перестановка $\sigma$ есть цикл $(12 \, \dots q)$, а на рис. \ref{Fig:REps}c перестановки $\rho_1, \rho_2, \rho_3$ задаются алгоритмом Ведюшкиной-Харчевой \cite{VedKha} для выбранного седлового атома $V$ (как со звездочками, так и без них, см. \cite{BolFom}).
	
	\subsection{Реализация числового инварианта расслоений Зейферта}
	
	В работе \cite{VedKib20} В.~В.~Ведюшкиной и В.~А.~Кибкало была подробно описана реализация биллиардными книжками $\Omega_k, \Omega_{k, k-s}$ для целых $0 \le s \le k$ произвольного значения метки $n = k$ на семье молекулы некоторого специального, зависящего от $k$ вида.
	
	В работе \cite{DAN20} также было описано семейство $\Omega_k'$ столов, для которых инварианты Фоменко-Цишанга при всех $k$ отличаются \textit{только} в метке $n = k$ одной семьи. Эти, в частности, удалось показать, что \textit{сложность} и \textit{валентность} семьи не являются (сами по себе) препятствиями к реализации семьи с меткой $n$ на ней с помощью биллиарда. Напомним, что сложностью семьи называют количество седловых особых окружностей у ее атомов, а валентностью атома (или графа из атомов) --- количество внешних ребер данного подграфа (в нашем случае --- количество граничных торов семьи).
	
	В настоящей работе мы изложим и наглядно проиллюстрируем идею построения биллиардных столов из работ \cite{VedKib20} и \cite{DAN20}.
	
	Напомним, что исследованных ранее интегрируемых биллиардах ненулевая метка $n$ возникала в двух основных случаях: когда область стола содержит фокусы семейства квадрик (например, биллиард в области $A_2$, ограниченной эллипсом), или когда слоение биллиарда имеет атомы со звездочками, а неплоский биллиардный стол --- коническую точку, т.~е. склейку пары столов по паре их общих дуг границы, образующих угол $\pi/2$.
	
	\begin{remark}
		Здесь и далее символом $\lambda$ обозначается как числовое значение биллиардного интеграла $\Lambda$ \eqref{eq_integral}, так и однозначно определенный цикл на граничном торе 3-атома.
	\end{remark}
	
\begin{figure}[!htb]
		\centering{\includegraphics[width=\linewidth]{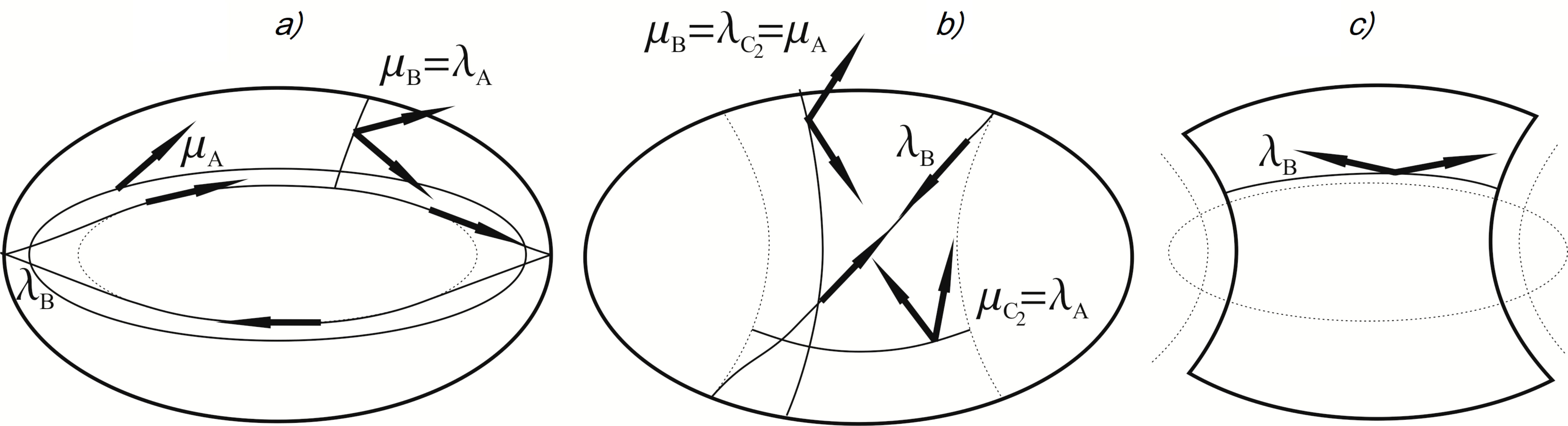}}
		\caption{Проекция циклов допустимых базисов на эллипс для биллиардов $A_2$ (рис. (а) и (b)), $A_0$ (рис. (c)).}
		\label{Fig:circles2}
	\end{figure}

	\begin{example}[Биллиард $A_2$ внутри эллипса] Для этой системы изобразим проекцию циклов $\lambda, \mu$ допустимых базисов (лежащих на регулярных граничных торах нужного атома) на биллиардный стол. Молекула имеет вид $2(A-) = B - A$, и особые слои атомов $2A, B, A$ лежат на уровнях $\lambda \in \{0, b ,a\}$ интеграла $\Lambda$.
		
		Циклы допустимых базисов на уровнях интеграла $\lambda \in \{\varepsilon, b - \varepsilon\}$ (cм. рис. \ref{Fig:circles2}a) и на уровнях  $\lambda \in \{b + \varepsilon, a - \varepsilon\}$ (см. рис. \ref{Fig:circles2}b) обозначим $(\lambda_A^{-}, \mu_A^{-}), (\lambda_B^{-}, \mu_B^{-})$ и $(\lambda_B^{+}, \mu_B^{+}), (\lambda_A^{+}, \mu_A^{+})$ соответственно. Отметим, что проекция цикла $\mu_A^{-}$ гомотопна граничному эллипсу. На уровне интеграла $b$ особая окружность атома $B$ проецируется на отрезок фокальной оси $Ox$, лежащий внутри эллипса.
		
		Направив стрелки на ребрах инварианта Фоменко-Цишанга от эллиптических атомов $A$ к седловому, имеем матрицы перехода $C_i^{\pm}$ на торах $\Lambda = b + \varepsilon$
		\[
		C_i^{-} =  \begin{pmatrix}  \lambda_{B}^{-} \\ \mu_B^{-} \end{pmatrix}  =
		\begin{pmatrix} 2 & 1 \\ 1 & 0  \end{pmatrix} \begin{pmatrix} \lambda_{A}^{-} \\ \mu_A^{-} \end{pmatrix}  \qquad
		C_i^{+} =  \begin{pmatrix}  \lambda_{B}^{+} \\ \mu_B^{+} \end{pmatrix}  =  \begin{pmatrix} -1 & 1 \\ 0 & 1
		\end{pmatrix}  \begin{pmatrix} \lambda_{A}^{+} \\ \mu_B^{+} \end{pmatrix}. \]
		
		Метка $n$ определяется как следующая сумма по всем ребрам семьи (см.~\cite{BolFom}):
		\[n = \sum_i [\cfrac{-\delta_i}{\beta_i}] + \sum_j [\cfrac{\alpha_j}{\beta_j}] + \sum_k [\cfrac{-\gamma_k}{\alpha_k}], \qquad C_s = \begin{pmatrix}\alpha_s & \beta_s \\ \gamma_s & \delta_s \end{pmatrix}. \]
		Здесь ребра с индексами $s = i, j, k$ --- входящие, выходящие и внутренние ребра семьи (в смысле поставленных выше стрелок) c матрицами перехода $C_i$ на них (выражающих базис 3-атома в конце ребра через базис 3-атома в начале ребра).
	\end{example}
	
	\begin{remark}
		Вклад в метку $n$ атома-семьи $B$ дает только верхнее ребро ($b< \lambda < a$) молекулы, и он равен $-1$.
	\end{remark}
	
	Ключевым является следующее построение. Пусть ребро молекулы со значением $b < \Lambda$, входящее атом фокального слоя, начинается не в максимальном атоме $A$, а в некотором седловом атоме $V$ на уровне $b < \lambda_1 < a$. Тогда (с точностью до знаков) имеем
	\[\lambda_V = \pm \mu_A^{+}, \quad \mu_V = \pm \lambda_A, \qquad C_i = \begin{pmatrix} -1 & 1 \\ 0 & 1
	\end{pmatrix} \begin{pmatrix} 0 & \pm 1 \\ \pm 1 & 0 \end{pmatrix}  = \begin{pmatrix} \pm 1 & \mp 1 \\ \pm 1 & 0 \end{pmatrix}.\]
	Тем самым, вклад ребра как входящего ребра в атом $B$ станет равен нулю, а как выходящего ребра из атома $V$ --- станет отличен от нуля и будет равен $\pm 1$.
	
	\begin{example}\label{Examp:121}
		Разрежем область $A_2$ на три части ветвями гиперболы $\lambda = \lambda_1 >b$. Обозначим левую и правую области $A_1$ с одним фокусом цифрами 1 и 4, а оставшуюся область без фокусов (тип которой обозначают $A_0$), лежащую между ветвями гиперболы, и ее дубликат --- цифрами 2, 3.
		
		Склеим книжку из листов $1, 2, 3, 4$ по перестановкам $\sigma = (1\,2\,3)$ на дуге левой ветви гиперболы, и $\rho = (2\,4\,3)$ --- на правой. Перестановки задаются простым правилом: на уровне $\Lambda = 0$ должна иметься пара минимальных особых окружностей (симметричных друг другу при отражении $y \to -y$), которые не меняют направления движения (т.~е. не отражаются назад от дуг склейки), и каждая из них проецируется на граничный эллипс стола $A_2$ биективно. Из $\sigma(1) = 2$ последовательно получим $\rho(2) = 4, \rho(4) = 3, \sigma(3) = 1$.
		
		Атом, возникающий вблизи уровня $\lambda = \lambda_1 >b$, имеет тип $C_2$. Он является прямым произведением 2-базы (2-атома $C_2$) на слой $S^1$ (лежащий в прообразе гиперболической дуги склейки). Эта вершина молекулы инцидентна двум атомам $B$ на уровне $\lambda = b$. Два эти атома (и целиком две компоненты связности слоения Лиувилля при $0 \le  \Lambda < b + \varepsilon$)  легко отождествить с такими же слоениями Лиувилля плоских биллиардов $A_2$ и $A_0$ (сохраняя интеграл $\Lambda$ и проекцию каждой точки каждого слоя на плоскость). Для биллиарда в области $A_0$ цикл $\lambda_B$ седлового атома изображен на рис.~\ref{Fig:circles2}c.
	\end{example}
	
	\begin{remark}
		Молекула имеет три одноатомных семьи. Метки $n$ на обоих атомах $B$ равны нулю. Метка $n$ на 3-атоме $C_2$ отлична от нуля и равна по модулю единице. Это есть вклад ребра, входящего в 3-атом $B$, проецирующийся на весь эллипс.
	\end{remark}
	
	В примере выше нам удалось ``запутать'' два стола $A_0$ и $A_2$ в единую книжку, так что слоения на уровне $0 \le \Lambda \le \lambda_1$ для $b < \lambda_1$ вошли в слоение результирующего биллиарда ``без изменений''.
	
\begin{figure}[!htb]
		\centering{\includegraphics[width=0.5\linewidth]{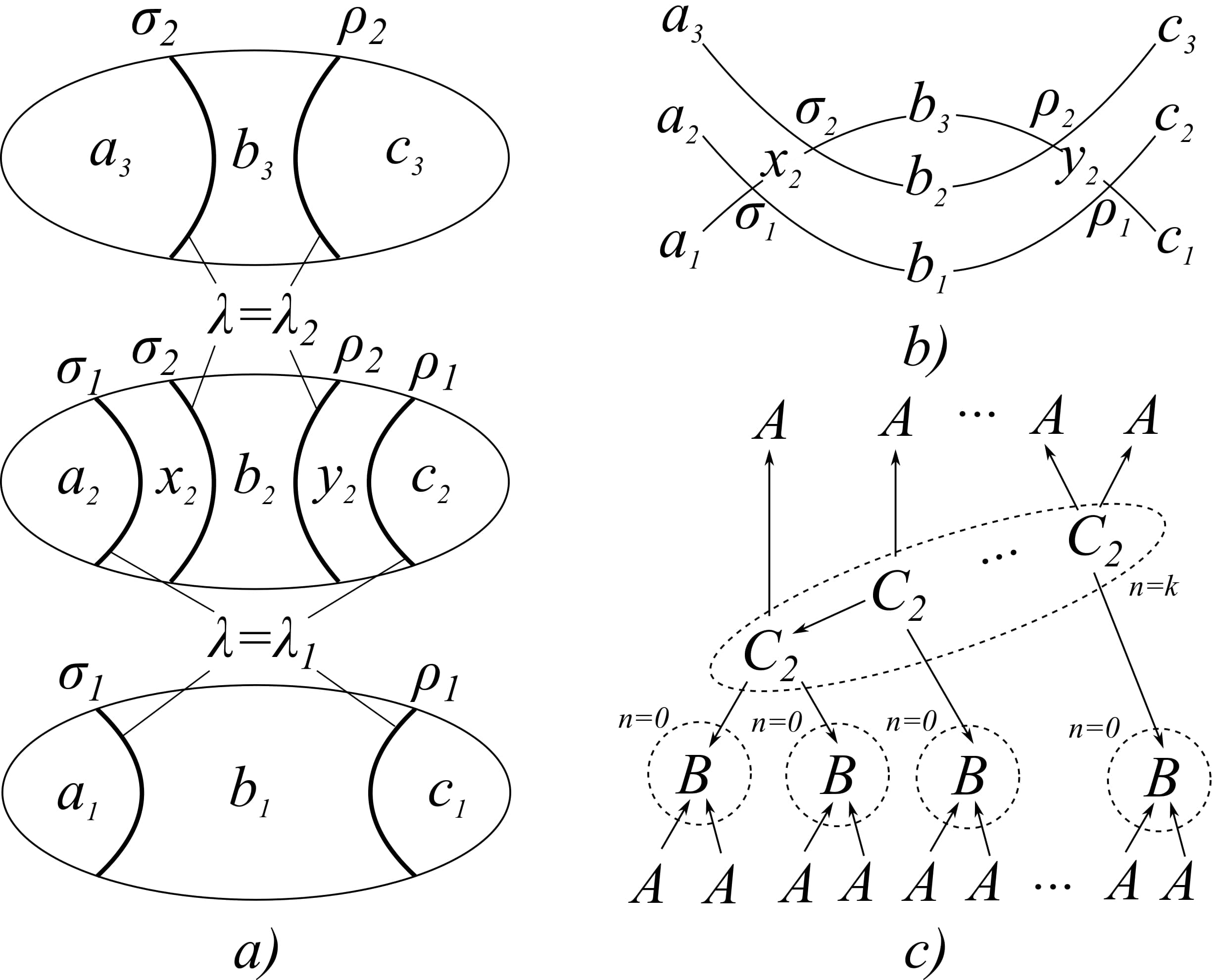}}
		\caption{Склейка биллиардной книжки $\Omega_3$ (а), ее схематичное изображение (b) и меченая молекула ее слоения Лиувилля (с).}
		\label{Fig:scheme_old}
	\end{figure}

	\begin{example}\label{Examp:222} Теперь разовьем эту конструкцию, и ``запутаем'' два стола-эллипса $A_2$. Аналогично разрежем столы-эллипсы по ветвям одной и той же софокусной им гиперболы $\lambda = \lambda_1$. Полученные куски столов $i = 1, 2$ обозначим $a_i^0, b_i^0, c_i^0$ (они совпадают с листами $a_1, b_1, c_1$  на рис. \ref{Fig:scheme_old}).
		
		По-прежнему потребуем, чтобы поднятие проекции биллиардного стола на плоскость порождало поднятие слоения плоского биллиарда в эллипсе при $0 \le \Lambda < \lambda_1$ на две связных компоненты слоения нового биллиарда, которое было бы послойным гомеоморфизмом (и лиувиллевой эквивалетностью).
		
		Тогда по известному началу кода траектории: $a_1, b_1, c_1$ или $a_2, b_2, c_2$ (кодом будет порядок столов, по которым проходит траектория), код траектории продолжается однозначно: $b_2, a_1$ или $b_1, a_2$. Каждая траектория должна вернуться на стол $a_i$ с тем же номером, с которого стартовала, причем ровно за один оборот. Кроме того, стол должен быть связен как клеточный комплекс.
		
		Полученные перестановки $\sigma$ и $\rho$ на левой и правой ветвях гиперболы имеют вид $\sigma = (a_2 \, b_2 \, a_1 \, b_1)$ и $\rho = (b_1 \, c_1 \, b_2 \, c_2)$. Грубая молекула не отличается от построенной в примере \ref{Examp:121}, а метка $n$ на семье $C_2$ равна $\pm 2$ из соображений симметрии.
	\end{example}
	
	\begin{remark}
		Здесь и далее листы $a_i, a_i^0$ лежат в полуплоскости $x <0$ (как и лист 1 из примера \ref{Examp:121}), а листы $c_i, c_i^0$ --- в полуплоскости $x >0$ (как и лист 4 из примера \ref{Examp:121}), и содержат по одному фокусу семейства квадрик (т.~е. имеют тип $A_1$). Листы $b_i$ и $b_i^0$ лежат между ветвями некоторой гиперболы (для всех $b_i^0$ гипербола имеет один и тот же параметр в семействе квадрик, софокусных с эллипсом границы).
	\end{remark}

	\begin{example}\label{Examp:kkk} Применим информацию их примера \ref{Examp:222} к случаю $k$ столов-эллипсов, разбитых на области $a_i^0, b_i^0, c_i^0$. По-прежнему потребуем связность стола и прохождение траекторией биллиарда столов $a_i^0, b_i^0, c_i^0$ именно в таком порядке. Точнее говоря, номеру $i$ сопоставлена пара симметричных относительно оси $Ox$ траекторий уровня $\Lambda  = 0$ и класс близких к ним при $0 < \Lambda < \varepsilon$, которые проходят столы комплекса в одном и том же порядке. Обе такие траектории должны замкнуться после прохождения по столу $b_{i+1}^0$ для номера $1 \le i < k$ или по столу $b_1$ для номера $k$. Тогда перестановки $\sigma, \rho$ имеют вид
		\[\sigma = (a_k \, b_k \, ... \, a_2 \, b_2 \, a_1 \, b_1 ) \qquad   \rho = (b_1 \, c_1 \, b_2 \, c_2 \,  ... \, b_k \, c_k).\]
		
		Грубая молекула имеет $k$ семей, состоящих из одного атома $B$, и одну семью типа прямого произведения, на 2-базе которого $k$ окружностей перестраиваются в $k$ окружностей через 2 критические точки. Слоение на базе этого 3-атома является неморсовским 2-атомом, т.~е. содержит две особые точки типа мульти-седла с $k$ пересекающими прямыми.
		
		Поскольку метка $n$ также имеет смысл класса Эйлера (т.~е. определяет препятствие к продолжению сечения с выбранных циклов на граничных торах семьи внутрь нее как многообразия Зейферта), то она не зависит от выбора слоения Лиувилля внутри семьи. В частности, метка $n$ сохранится при возмущения атома в семье с его распадом на молекулу из боттовских атомов. Тем самым на таком столе из $k$ областей $A_2$ уже реализована метка $n = k$, пусть и слоением Лиувилля с неботтовским атомом.
	\end{example}
	
	Последним шагом явно опишем устройство стола $\Omega_k$ после возмущения. На рисунке \ref{Fig:scheme_old}a изображена склейка  книжки $\Omega_3$ из частей трех столов-эллипсов. На рис. \ref{Fig:scheme_old}b изображена схема ее сечения $x = 0$, т.~е. ``вид сбоку''. Возмущение неботтовского атома дает граф-дерево из $k-1$ атомов $C_2$, соединенных в цепочку.
	
	Требуемое разрезание $k$ столов $A_2$ описано в таблице \ref{Tab:book_sheet}. Cтол с номером $1 < i < k$ будем разбивать на 5 частей $a_i, x_i, b_i, y_i, c_i$, границами которых являются ветви гипербол $\lambda_{i-1} < \lambda_{i}$. Здесь $b < \lambda_1 < \dots < \lambda_{k-1} < a$. Столы с номерами $1$ и $k$ разобьем ветвями лишь одной гиперболы, с параметром $\lambda_1$ и $\lambda_{k-1}$ соответственно. Для столов $\Omega_k$ из работы \cite{VedKib20} перестановки $\sigma_i$ и $\rho_i$ на ветвях гиперболы $\lambda_{i}$ указаны в таблице \ref{Tab:permutations}.

	\begin{table}
		\begin{center}
			\begin{tabular}[t]{|c||c|c|c|c|c|c|c|c|c|c|c|c|c|}
				\hline
				Область & Тип & Номер стола & Граница & $Oxy$  \\
				\hline
				$a_i$ & $A_1$  & $S_i, 1 \le i \le k$ & $\lambda_{1}$ при $i = 1$\,; $\lambda_{i-1}$ при $2 \le i \le k$ & $x < 0$   \\
				\hline
				$x_i$ & $A_0$  & $S_i, 2 \le i \le k-1$ & $\lambda_{i-1}$\, и $\lambda_{i}$ при $2 \le i \le k-1$ & $x < 0$   \\
				\hline
				$b_i$ & $A_0$  & $S_i, 1 \le i \le k$ & $\lambda_{1}$ при $i = 1$\,; $\lambda_{i-1}$ при $2 \le i \le k$ & $Oy \subset b_i$ \\
				\hline
				$y_i$ & $A_0$  & $S_i, 2 \le i \le k-1$ & $\lambda_{i-1}$\, и $\lambda_{i}$ при $2 \le i \le k-1$ & $x > 0$   \\
				\hline
				$c_i$ & $A_1$  & $S_i, 1 \le i \le k$ & $\lambda_{1}$ при $i = 1$\,; $\lambda_{i-1}$: при $2 \le i \le k$ & $x > 0$   \\
				\hline
			\end{tabular}
			\label{Tab:book_sheet}
			\caption{Обозначения листов биллиардных столов.}
		\end{center}
	\end{table}
	
	\begin{theorem}[В.~В.~Ведюшкина, В.~А.~Кибкало]
		По любому $k \in \mathbb{Z}$ алгоритмически строится биллиард $\Omega_k$ со слоением Лиувилля на неособой изоэнергетической поверхности, содержащим некоторую семью с меткой $n = k$. Его инвариант Фоменко-Цишанга изображен на рис.~\ref{Fig:scheme_old}c.
	\end{theorem}
	
	\begin{table}[ht]
		\begin{center}
			\begin{tabular}[t]{|c||c|c|c|c|c|c|c|c|c|c|c|c|c|}
				\hline
				Гипербола & $\rho_i$ & $\sigma_i$ \\
				\hline
				$ i = 1$ & $(b_1, c_1, y_2, c_2)$ & $(a_2, x_2, a_1, b_1)$   \\
				\hline
				$ 2 \le i \le k-1$ & $(b_{i}, y_{i}, y_{i+1}, c_{i+1})$ & $(a_{i+1}, x_{i+1}, x_i, b_i)$   \\
				\hline
				$ i = k-1$ & $(b_{k-1}, y_{k-1}, b_{k}, c_{k})$ & $(a_{k}, b_{k}, x_{k-1}, b_{k-1})$   \\
				\hline
			\end{tabular}
			\label{Tab:permutations}
			\caption{Перестановки на корешках склейки столов $\Omega_k$.}
		\end{center}
	\end{table}
	
	Модифицировав столы $\Omega_k$, можно реализовать инвариант Фоменко-Цишанга, на котором указанная метка не обязательно равна $k$, но не превосходит по модулю $k$. Для этого удалим из $s$ столов $A_2$ области $a_i, c_i$, и ``пропустим'' их при записи перестановок $\sigma_i, \rho_i$. Так строятся столы из серии $\Omega_{k, k-s}$. Отметим, что $\Omega_{k, k} = \Omega_k$.
	
	\begin{theorem}[В.~В.~Ведюшкина, В.~А.~Кибкало]
		При любых целых $k, s$, где $0 \le s \le k$, инвариант Фоменко-Цишанга cлоения Лиувилля книжки $\Omega_{k, k-s}$ совпадает с инвариантом $\Omega_{k}$ за исключением метки $n = k-s$ на семье из $k-1$ атома $C_2$, имеющей валентность $2k$. Остальные метки и атомы не зависят от $s$.
	\end{theorem}
	
	\begin{remark}
		В примерах \ref{Examp:121} и \ref{Examp:222} были описаны столы $\Omega_{2, 1}$ и $\Omega_{2, 2} = \Omega_2$ и объяснено отличие в значениях метки $n$. Общий случай получается по аналогии.
	\end{remark}

\begin{figure}[!htb]
		\centering{\includegraphics[width=0.85\linewidth]{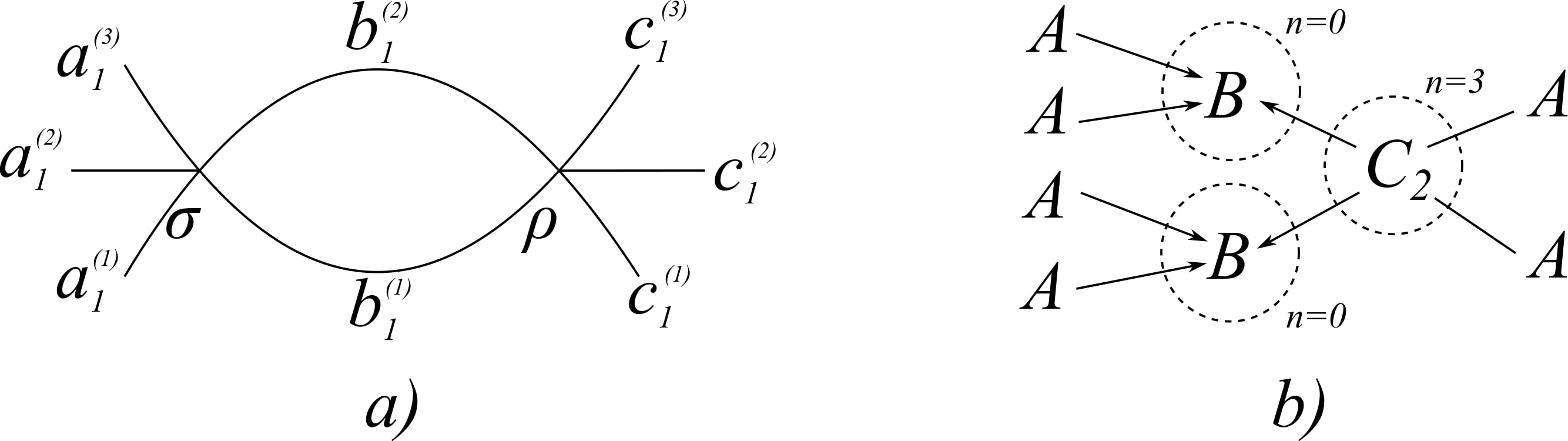}}
		\caption{Схематичное изображение биллиардной книжки $\Omega'_3$ (a) и меченая молекула ее слоения Лиувилля (b).}
		\label{Fig:n_new}
	\end{figure}

	Теперь опишем конструкцию книжек $\Omega_k'$, введенных в работе \cite{DAN20}. Возьмем $k$ областей $a_1^0$, $k$ областей $c_1^0$ и две области $b_1^0$. Обозначим их так: $a^{i}_1, b^{1}_1, b^{2}_1, c^{i}_1$ для $i = 1, \dots, k$. Выберем такие перестановки, чтобы траектория на уровне $\Lambda < b$ сначала двигалась над левой областью $A_1$ плоскости, проходя по всем $a_1^i$ от $a_1^1$ до $a_1^k$ (отражаясь несколько раз от границы по перестановке $\sigma_i$), затем переходила на стол $b_1$, далее двигалась аналогично по столам $c_1^i$, и в конце возвращалась через стол $b_2$. При $k =3$ схема книжки и ее инвариант Фоменко-Цишинга изображены на рис. \ref{Fig:n_new}a и \ref{Fig:n_new}b.
	
	\begin{theorem}
		Для каждого неотрицательного $k \in \mathbb{Z}$ книжки $\Omega_k'$, склеенные из листов $a^{i}_1, b^{1}_1, b^{2}_1, c^{i}_1, 1 \le i \le k$ по перестановкам \[\sigma = (a^1_1 \, a^2_1 \, \dots \, a^k_1 \, b^1_1 \, b^2_1), \quad \rho = (b^1_1 \, c^1_1 \, c^2_1 \, \dots \, c^k_1 \, b^2_1),\] реализуют слоения Лиувилля, чей инвариант Фоменко-Цишанга изображен на рисунке \ref{Fig:scheme_old}c. Метка $n = k$ на семье из одного атома $C_2$. Остальные элементы инварианта не зависят от $k$.
	\end{theorem}
	
	Как следствие, метка на семье инварианта $\Omega_k'$ (в отличие от серий $\Omega_{k, k-s}$)  не ограничена валентностью или сложностью семьи --- любое возможное значение $k$ метки $n$ реализуется на достаточно простой семье, одной и той же для всех $k$. Будет интересно проверить, какие именно семьи из одного или нескольких атомов можно получить на месте $C_2$ при развитии данной конструкции.
	
	\section{Свойста слоений биллиардных систем}\label{s:4:exotic}
	
	Данный раздел мотивирован следующим общим вопросом: какие свойства и утверждения о топологии и динамике гладких и вещественно-аналитических систем остаются верными (или имеют аналоги) в биллиардных системах. Например, обязательно ли выполнен аналог теоремы Лиувилля для слоения Лиувилля биллиардной книжки, перестановки которой коммутируют?
	
	Дело в том, что фазовое пространство биллиардной системы получается отождествлением точек на границе некоторого множества. Гладкая и симплектическая структура при этом не определялась в прообразе границы и склеек стола. Там склейка гарантирует лишь непрерывность функций и слоения Лиувилля. На вопрос о гамильтоновом сглаживании (т.~е. введения гладкой и симплектической структур вблизи прообраза границы) частичный положительный ответ получен В.Ф. Лазуткиным \cite{Lazutkin} в случае склейки стола по выпуклым дугам. Отметим также важные продвижения в вопросе изучения интегрируемых геодезических потоков и биллиардов на склеенных пространствах, полученные Е.~А.~Кудрявцевой в работе \cite{KudrBill}.
	
	\subsection{Топологическая устойчивость}\label{s:4:1:flow}
	
	Определение инварианта Фо\-мен\-ко-Ци\-шан\-га слоения Лиувилля на 3-подмногообразии $Q^3_h$ уровня энергии $H = h$ требовало ряда ограничений на ИГС. Опыт исследования реальных систем физики и механики показал, что они являются вполне естественными и выполняясь практически всегда. Так, для критических окружностей в $Q^3_h$ должно выполняться введенное А.Т.Фоменко свойство Морса--Ботта. Затем, система должна быть нерезонансной, т.~е. почти на всех слоях частоты интегралов $H$ и $F$ несоизмеримы.
	
	Кроме того, ИГС $\mathfrak{F} = (H, F)$  на симплектическом многообразии $(M^4, \omega)$ должна быть \textit{топологически устойчивой} на уровне $Q^3_h$. По определению это значит, что при любом малом изменении $\varepsilon$ значения энергии $H$ структура слоения Лиувилля на $Q^3_{h+\varepsilon}$ не меняется.
	
	Тем самым, для доказательства возможности ввести инвариант Фоменко--Цишанга для некой биллиардной книжки требуется проверить топологическую устойчивость системы, а не только послойную гомеоморфность слоев торам, а особенностей --- 3-атомам.
	
	Для вещественно-аналитических ИГС топологическая устойчивость на $Q^3_h$ требует совпадения ориентации (гамильтоновым полем $\mathrm{sgrad}\,H$) всех критических окружностей, лежащих в одном седловом 3-атоме. Обратно, наличие 3-атома с разной ориентацией на таких окружностях гарантирует топологическую неустойчивость указанной системы. Покажем, что в биллиардах без потенциала такая неустойчивость невозможна:
	
	\begin{assertion}\label{Ass:1}
		Пусть слоение Лиувилля биллиарда c гамильтонианом $H = |\vec{v}^2|$ на некоторой биллиардной книжке при $H = h >0$ состоит из семейств регулярных торов и конечного числа особых слов, окрестности которых гомеоморфны 3-атомам. Тогда ориентация на особых окружностях седлового атома (отличного от $B, A^{*}$) одинакова.
	\end{assertion}
	
	\begin{proof}
		1. Особые окружности седловых 3-атомов биллиардов без потенциала проецируются либо на фокальную прямую, либо на дугу невырожденной квадрики --- эллипса или гиперболы. Прообраз такой дуги в особом слое --- одна или несколько точек, различаемых напавлением вектора скорости.
		
		Пусть это дуга эллипса $0 < \lambda_1 <b$. Отобразим симметрично относительно $Ox$ в верхнюю полуплоскость ту часть биллиарднного стола, которая проецируется на область $y <0$, добавив при необходимости транспозиции на новых ребрах склейки. Теперь все особые окружности уровня $\lambda_1$ проецируются на одну связную дугу в проекции стола.
		
		Рассмотрим прообраз гиперболы, пересекающей эту дугу во внутренней точке. Компоненты связности ее прообраза (в 3-атоме) гомеоморфны либо 2-базе этого 3-атома, либо его дублю (из вида слоения). При любом движении по прообразу этой гиперболы (в 3-атоме) сохраняется направление вектора скорости ``вправо'' или ``влево'', т.~е. биллиардное движение одинаково ориентирует особые окружности 3-атомов.
		
		2. На гиперболических уровнях $b< \lambda_2 < a$ применим отражение относительно $Oy$. Если фокус эллипса не входит в проецию книжки, то конструкция по-прежнему работает. В ином случае отобразим биллиардный стол в одну полуплоскость относительно $Oy$. Проекция стола содержит точки отрезка между фокусами. Для прообраза гиперболы, проходящей через такую точку, используем аналог вышеописанной конструкции.
	\end{proof}
	
	\begin{corollary}
		1. Пусть слоение Лиувилля некоторой биллиардной книжки послойно гомеоморфно слоению Лиувилля некоторой ИГС. Тогда для такого биллиарда определен инвариант Фоменко-Цишанга.
		
		2. ИГС на $Q^3_h$ с разной ориентацией окружностей одного из своих 3-атомов не моделируется в классе интегрируемых биллиардов без потенциала.
	\end{corollary}
	
	Согласно комментарию Е.А.Кудрявцевой, для достаточно широкого класса биллиардных книжек указанный вид слоения Лиувилля (гомеоморфность неособых слоям торам и послойная гомеоморфность окрестностей особых слоев боттовским 3-атомам) и одинаковая ориентация особых окружностей будут следовать из теорем, доказанных для гладких и вещественно-аналитических систем (приведенных, например, в \cite{BolFom}) благодаря результатам В.Ф.Лазуткина \cite{} о существовании гамильтонова сглаживания. В этом случае важно, чтобы траектории на особом слое и близких слоях не касались граничных невыпуклых дуг склейки.
	
	Отметим, что интегрируемое возмущение такой системы, при котором окружности с разной ориентацией окажутся на разных уровнях интеграла, будет топологически устойчивой системой. Слоение Лиувилля такой возмущенной системы будет иметь инвариант Фоменко--Цишанга и включать в себя семью из нескольких седловых атомов с меткой $\varepsilon = -1$ на ребре между ними. Согласно утверждению \ref{Ass:Nomer2}, в классе круговых биллиардов такие системы реализуются.
	
	Примером является стол, склеенный из трех областей: двух одинаковых круговых колец $\alpha, \beta$ с радиусами $0 < r_1 < r_2$ и диска $\gamma$ радиуса $r_1$ с перестановкой $(\alpha\, \beta\, \gamma)$. Критические окружности двух атомов $B$ проецируются на окружность склейки и отличаются ориентацией. Любая из двух деформаций линии склейки (при достижении центра семейства окружностей или внешней граничной окружности) разрушает оба атома $B$.
	
	Тем самым, указанный стол не допускают деформацию, при которой седловые атомы семьи преобразуются в единый атом (и их особые окружности попадают на один уровень интеграла). Отметим, что такого же эффекта в \textit{софокусных} биллиардах не наблюдалось. Сформулируем возникший вопрос, в дополнение к локальной гипотезе $ \mathbf{C}$:
	
	\begin{hypothesis}
		В классе биллиардных книжек, склеенных из софокусных биллиадов, не встречается числовой инвариант $r = \infty, \varepsilon = -1$ на ребре между седловыми атомами.
	\end{hypothesis}
	
	\subsection{Расщепляющиеся 3-атомы в биллиарде с потенциалом Гука}\label{s:4:2:split}
	Рассмотрим особый слой $L$, не имеющий точек ранга 0, но содержащий несколько седловых окружностей $S^1_i$ точек ранга 1. Пусть вблизи $L$ критическое множество $X$ имеет вид произведения отрезка $I$ на $\sqcup S^1_i$. Тогда бифуркационная диаграмма $\Sigma$ --- образ $X$ при отображении момента --- локально имеет вид одной кривой, или нескольких пересекающихся кривых. Во последнем случае особый слой $L$ распадается на несколько особых слоев при небольшом изменении $h$. Такую особенность $L$ называют \textit{расщепляющейся}. Нам будет достаточно данного простого примера.
	
	Покажем, что в слоении Лиувилля биллиардов с потенциалом встречаются расщепляющиеся особенности ранга 1. Напомним, что добавление потенциала Гука к системе биллиарда в софокусных квадриках сохраняет интегрируемость, что следует из работы В.~В.~Козлова \cite{Kozlov}. Топология слоений Лиувилля плоских биллиардов с потенциалом Гука изучалась в работах И.~Ф.~Кобцева (биллиард внутри элллипса \cite{HookeEllipse}), С.~Е.~Пустовойтовым (биллиард в кольце между эллипсами \cite{HookeRing}).
	
	В эллиптических координатах $\lambda_1, \lambda_2$ и сопряженных импульсах $\mu_1, \mu_2$ энергия $H$ и дополнительный интеграл $F$ системы имеют вид
	\[H = \cfrac{2(\lambda_1 + a)(\lambda_1 + b)}{\lambda_2 - \lambda_1}\,\mu_1^2 + \cfrac{2(\lambda_2 + a)(\lambda_2 + b)}{\lambda_1 - \lambda_2} \mu_2^2 + \cfrac{k}{2} (\lambda_1 + \lambda_2),\]
	\[F = 2(\lambda_1 + a)(\lambda_1 + b) \mu_1^2 + \cfrac{k}{2}\lambda_1^2 - \lambda_1 H.\]
	Знак числа $k$ определяет притягивающий ($k>0$) или отталкивающий ($k<0$) потенциал.
	
\begin{figure}[!htb]
		\includegraphics[width=\linewidth]{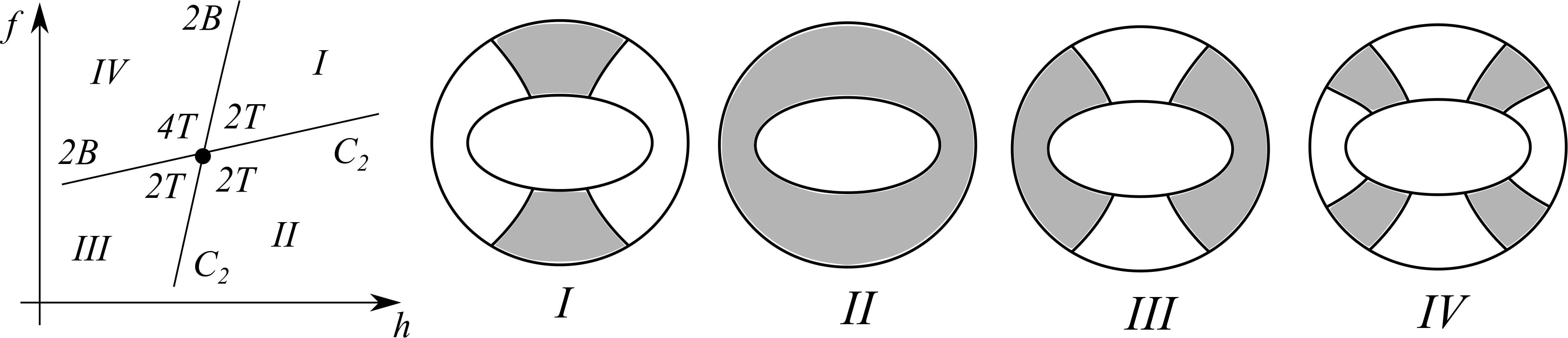}
		\caption{Локальная бифуркационная диаграмма и проекции уровня $H = h, F = f$ на плоскость $Oxy$ для камер I-IV.}
		\label{Fig:Hook}
	\end{figure}

	При наличии потенциала Гука топология слоения Лиувилля система, ограниченная на $Q^3_h$, будет зависеть от выбора $h$. Бифуркационая диаграмма $\Sigma$ зависит от знака $k$. Нам интересен случай $k< 0$. При $k <0$ одна из особых точек бифуркационной диаграммы $\Sigma$ системы биллиарда внутри эллипса имеет круговую молекулу (инвариант слоения на 3-границе, см. \cite{BolFom}) как у 4-особенности седло-седло $(B \times C_2) \slash \mathbb{Z}_2$ \cite{HookeEllipse}.
	Переход к системе в эллиптическом кольце не меняет кривых бифуркационной диаграммы при $k <0$, но изменяет типы 3-атомов, им соответствующие. Для четырех камер $I, II, III, IV$ плоскости $\mathbb{R}^2(h,f)$ покажем на рис. \ref{Fig:Hook}, как устроена проекция поверхности уровня $(h, f)$ пары функций $H, F$ на плоскость $Oxy$.
	
	Круговая молекула упомянутой выше особой точки $\Sigma$ тоже изменится, и будет отличаться от круговых молекул особенностей ранга 0 типа седло-седло. Рассмотрим прообразы малых отрезков, проходящих через особую точку $\Sigma$ трансверсально обеим пересекающимся в ней кривым. Такие отрезки принадлежат одному из двух классов, что определяется выбором пары вертикальных углов в точке пересечения.
	
	\begin{assertion}
		Данная особенность послойно гомеоморфна расщепляющейся особенности ИГС. Слоение Лиувилля в прообразе малого вертикального отрезка имеет тип 3-атома $P^4$ (между камерами II и IV), а в прообразе кривых из другого класса --- тип 3-атома $K_2$ (между камерами III и I).
	\end{assertion}
	
\begin{figure}[!htb]
		\includegraphics[width=0.9\linewidth]{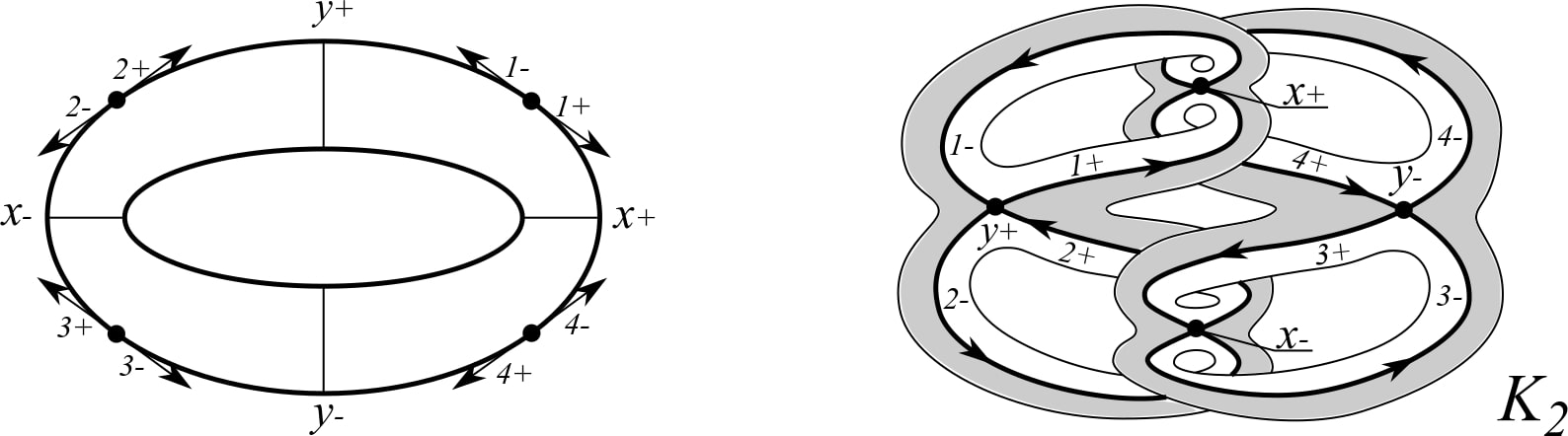}
		\caption{Дуги и точки эллипса соответствуют лентам $1\pm, 2\pm, 3\pm, 4\pm$ и седлам $x\pm, y\pm$ атома $K_2$ при переходе из камеры III в камеру I.}
		\label{Fig:K2}
	\end{figure}

	\begin{proof}
		1. В прообразе точки пересечения дуг $\Sigma$ лежит связный двумерный слой с четырьмя особыми окружностями. Они гомологичны друг другу на слое, а их проекция $\pi$ на плоскость $Oxy$ переводит их в четыре отрезка осей $Ox, Oy$, лежащие в области $\Omega$.
		
		Поскольку все седловые атомы здесь имеют тип прямого произведения, то 4-особенность есть произведение окружности $S^1$ на трехмерную базу. Тем самым, особенность послойно гомеоморфна расщепляющейся особенности ИГС, а ограничение на $Q^3_{\gamma} = \mathfrak{F}^{-1}(\gamma)$ --- некоторому атому ИГС с 1 степенью свободы.
		
		2. На рис.\ref{Fig:K2}a оснастим касательными векторами по и против часовой стрелки точки внешней эллиптической границы закрашенных областей. Отождествив пары этих векторов в угловых точках закрашенной области, получим кольца атома двух цветов (белые и черные). Белым цветом обозначим кольца той камеры, где значение $H$ меньше (в паре камер III и I это будет камера III).
		
		Концы полуосей эллипса (как точки проекции критических окружностей 3-атома на стол) соответствуют седловым особым точкам атома, т.~е. его крестам. Ленты атома здесь соответствуют четверти эллипса между двумя концами его полуосей (обозначим цифрами 1, 2, 3, 4) с указанием направления обхода: знаки ``$+$'' и ``$-$'' для направлений по и против часовой стрелки соответственно.
		
		3. Атом в прообразе вертикального отрезка имеет симметрию $\mathbb{Z}_4$ на вершинах и перестраивает два кольца в четыре. Тогда его тип --- $P_4$. Другой тип кривых соответствует атому $K_2$, см. рис. \ref{Fig:K2}.
		
	\end{proof}

\end{document}